\documentclass[11pt]{article}

\usepackage{amssymb}
\usepackage{amsmath}
\usepackage{amstext}
\usepackage{amsthm}
\usepackage{fullpage}
\usepackage[dvips]{graphicx}
\usepackage{subfigure}              % edo

\newcommand{\C}{\mathbb C}
\newcommand{\R}{\mathbb R}

\newcommand{\Z}{\mathbb Z}
\newcommand{\E}{{\mathcal E}}
\newcommand{\HH}{{\mathcal H}}

\newcommand{\eps}{\varepsilon}

\newcommand{\set}[1]{\left\{#1\right\}}

\newtheorem{theorem}{Theorem}[section]

\newtheorem{prop}[theorem]{Proposition}
\newtheorem{cor}[theorem]{Corollary}
\theoremstyle{definition}
\newtheorem{definition}[theorem]{Definition}

\theoremstyle{remark}

\numberwithin{equation}{section}

\begin{document}

\title{Hylomorphic solitons  in the  nonlinear Klein-Gordon equation }
\author{J. Bellazzini  \thanks{Dipartimento di Matematica Applicata,
Universit\`a degli Studi di Pisa, Via F. Buonarroti 1/c, Pisa, ITALY.
e-mail: \texttt{j.bellazzini@ing.unipi.it, benci@dma.unipi.it, bonanno@mail.dm.unipi.it}}
\and V. Benci\addtocounter{footnote}{-1}\footnotemark \and C. Bonanno\addtocounter{footnote}{-1}\footnotemark \and E. Sinibaldi \thanks{Scuola Superiore S. Anna, CRIM-Lab, Viale R. Piaggio 34, Pontedera (Pisa), ITALY. e-mail: \texttt{edoardo.sinibaldi@crim.sssup.it}}}

\date{}
\maketitle
\abstract{Roughly speaking a solitary wave is a solution of a field equation whose
energy travels as a localised packet and which preserves this localisation
in time. A soliton is a solitary wave which exhibits some strong form of
stability so that it has a particle-like behaviour. In this paper we show a new mechanism
which might produce solitary waves and solitons for a large class of equations, such as the nonlinear Klein-Gordon equation. We show that the existence of these kind of solitons, that we have called \emph{hylomorphic} solitons, depends on a suitable energy/charge ratio. We show a variational method that allows to prove the existence
of hylomorphic solitons and that turns out to be very useful for numerical applications. Moreover we introduce some classes of nonlinearities which admit hylomorphic solitons of different shapes and with different relations between charge, energy and frequency.}
\tableofcontents

\section{Introduction}

Roughly speaking a solitary wave is a solution of a field equation whose
energy travels as a localised packet and which preserves this localisation
in time.

A \textit{soliton} is a solitary wave which exhibits some strong form of
stability so that it has a particle-like behaviour.

Today, we know (at least) three mechanism which might produce solitary
waves and solitons:

\begin{itemize}
\item  Complete integrability; e.g. Korteweg-de Vries equation
\begin{equation}
u_{t}+u_{xxx}+6uu_{x}=0  \tag{KdV}  \label{KdV}
\end{equation}

\item  Topological constraints: e.g. Sine-Gordon equation
\begin{equation}
u_{tt}-u_{xx}+\sin u=0  \tag{SGE}  \label{SG}
\end{equation}

\item  Ratio energy/charge: e.g. the following nonlinear Klein-Gordon
equation
\begin{equation}
\psi _{tt}-\Delta \psi +\frac{\psi }{1+\left| \psi \right| }=0 \qquad \qquad \psi \in
\mathbb{C}  \tag{NKG}  \label{NKG}\\
\end{equation}
\end{itemize}

This paper is devoted to the third type of solitons which will be
called \emph{hylomorphic solitons}. This class of solitons that
are characterised by a suitable energy/charge ratio includes the
$Q$-balls, spherically symmetric solutions of \ref{NKG}, as well
as solitary waves which occur in the nonlinear Schr\"odinger
equation (see e.g. \cite{CL82}, \cite{BBGM}) and in gauge theories
(see e.g. \cite{bf}, \cite{befo}, \cite{befogranas}). We have chosen
the name \emph{hylomorphic}, which comes from the Greek words
``\textit{hyle}''=``\textit{matter}''=``\textit{set of
particles}'' and ``\textit{morphe}''=``\textit{form}'', with the
meaning of solitons ``giving a suitable form to condensed matter''
(see Section \ref{def-hylo-sol} for definitions and details).

The aims of this paper are the following:

\begin{itemize}

\item to give the definition of \emph{hylomorphic solitons} and to set this notion in the literature on non-topological solitons;

\item to describe a new variational approach that allows to prove the
existence of hylomorphic solitons for a large class of nonlinearities for
\ref{NKG}. This variational method turns out to be very useful for numerical
simulations;

\item to classify the nonlinearities which give hylomorphic
solitons: this classification is based on the different shapes of the
solitons and the different relations between charge, energy and frequency.
We obtain four classes of nonlinearities and we prove necessary conditions
for the nonlinear term to belong to a given class. Numerical simulations show
the different behaviour of these classes quantitatively.
\end{itemize}

\section{The abstract theory}

\subsection{An abstract definition of soliton} \label{be}

We consider dynamical systems with phase space $X$ described by one or more fields,
which mathematically are represented by a function
$$
\Psi :\R^{n} \rightarrow V, \qquad n\ge 2
$$
where $V$ is a vector space with norm $\left\| \cdot \right\| _{V}$ which is
called the internal parameters space. We assume the system to be deterministic, and denote by $U:\R \times X \rightarrow X$ the time evolution map, which is assumed to be defined for all $t\in \R$. The dynamical system is denoted by $(X,U)$. If $\Psi _{0}(x)\in X$ is the initial condition, the evolution of the system is described by
$$
\Psi \left( t,x\right) =U_{t}\Psi _{0}(x)
$$
We assume that $X\subset L^{2}(\R^{n},V),$ namely $\int \left\| \Psi
\left( x\right) \right\| _{V}^{2}dx<+\infty $ for every $\Psi \in X$. Moreover for states $\Psi \in X$ which satisfy $\int \left| x\right| \left\| \Psi \left( x\right) \right\|_{V}^{2}dx<+\infty$, it is possible to give the notion of \textit{barycenter} of the state as follows
\begin{equation} \label{barycenter}
q(t)=q_{\Psi }(t)=\frac{\int x \left\| \Psi \left(
t, x\right) \right\| _{V}^{2}dx}{\int \left\| \Psi \left( t, x\right) \right\|
_{V}^{2}dx}
\end{equation}
The term \emph{solitary wave} is usually used for solutions of field equations whose energy is localised and the localisation of the energy packet is preserved under the evolution. Using the notion of barycenter, we give a formal definition of solitary wave.

\begin{definition} \label{solw}
A state $\Psi \in X$ is called \emph{solitary wave} if for any $\varepsilon >0$ there exists a radius $R>0$ such that for all $t\in \R$
$$
\int \left\| U_{t}\Psi \left( x\right) \right\| _{V}^{2}dx-\int_{B_{R}( q_{\Psi }(t))}\left\| U_{t}\Psi \left( x\right) \right\| _{V}^{2}dx<\varepsilon
$$
where $B_{R}( q_{\Psi }(t))$ denotes the ball in $\R^{n}$ of radius $R$ and centre in $ q_{\Psi }(t)$.
\end{definition}

The \emph{solitons} are solitary waves characterised by some form of stability. To
define them at this level of abstractness, we need to recall some well known
notions in the theory of dynamical systems.

\begin{definition} \label{stable-invariant}
Let $X$ be a metric space and let $\left( X,U\right) $ be a dynamical
system. An invariant set $\Gamma \subset X$ is called \emph{stable}, if  for any $\varepsilon >0$ there exists $\delta >0$ such that if $d(\Psi,\Gamma )< \delta$ then
$d(U_{t}\Psi ,\Gamma )< \varepsilon$ for all $t\in \R$.
\end{definition}

\begin{definition} \label{dos}
A state $\Psi_{0}$ is called \emph{orbitally stable} if there exists a
finite dimensional manifold $\Gamma \subset X$ with $\Psi _{0}\in \Gamma$, such that $\Gamma$ is invariant and stable for the dynamical system $(X,U)$.
\end{definition}

The above definition needs some explanation. Since $\Gamma $ is invariant, $U_{t}\Psi _{0}$ is in $\Gamma $ for all $t$. Thus, since $\Gamma $ is finite
dimensional, the evolution of $\Psi _{0}$ is described by a finite number of
parameters$.$ Thus the dynamical system $\left( \Gamma ,U\right) $\ behaves
as a point in a finite dimensional phase space. By the stability of $\Gamma$, the evolution $U_{t}\Psi_{\eps}$ of a small perturbation $\Psi_{\eps}$ of $\Psi _{0}$ might become very different from $U_{t}\Psi _{0}$, but $U_{t}\Psi_{\eps}$ remains close to $\Gamma$. Thus, the perturbed system appears as a finite dimensional system with a small perturbation which depends on an infinite number of parameters.

\begin{definition} \label{ds}
A state $\Psi _{0}\in X$,$\Psi _{0}\neq 0,\;$ is called \emph{soliton} if it is a orbitally stable solitary wave.
\end{definition}

According to this definition a soliton is a state in which the mass is
``concentrated'' around the barycenter $q(t)\in \R^{n}$. In
general, $\dim \left( \Gamma \right) >n$ and hence, the ``state'' of a
soliton is described by $n$ parameters which define its position and other
parameters which define its ``internal states''.

\subsection{Definition of hylomorphic solitons} \label{def-hylo-sol}

In this section we will expose a general method to prove the existence of non-topological
solitons. This method leads in a
natural way to the definition of \textit{hylomorphic solitons}.

We make the following assumptions: (i) there are at least two integrals of motion, the energy $\E$ and the \emph{hylomorphic charge} $\HH$; (ii) the system is invariant for space translations. In some models the hylomorphic charge is just the ``usual'' charge and this fact justifies this name. We add the attribute \emph{hylomorphic} just to recall that $\HH$ might not be a charge (as in the case of the
nonlinear Schr\"odinger equation). Moreover, in many models in quantum field
theory, $\HH$ represents the expected number of particles, hence  for charged particles it is proportional to the electric charge.

Now, given the set
$$
{\mathcal M}_{\sigma} =\left\{ \Psi \in X\ :\ \HH(\Psi)=\sigma \right\}
$$
In order to prove the existence of solitary waves and solitons, we follow a method based on the following steps:

\begin{itemize}
\item  (S-1) prove that the energy $\E$ has a minimum $\Psi _{0}$ on ${\mathcal M}_{\sigma}$;

\item  (S-2) prove that the barycenter $q_{\Psi}$ (see
(\ref{barycenter})) is well defined and that for any $\varepsilon
>0$ there exists $R>0$ such that
$$
\int \left\| \Psi \left( x\right) \right\|
_{V}^{2}dx-\int_{B_{R}(q_{\Psi })}\left\| \Psi \left(
x\right) \right\| _{V}^{2}dx<\varepsilon \qquad \forall\, \Psi \in
\Gamma
$$
where $\Gamma$ is the set of minimisers;

\item  (S-3) prove that $\Gamma $ is finite dimensional;

\item  (S-4) prove that $\Gamma $ is stable.
\end{itemize}

These steps will be explained in details for the nonlinear Klein-Gordon equation in Section \ref{qball}.

The integrals of motion are used in the definition of hylomorphic solitons. Setting
\begin{equation} \label{brutta}
m_{0}=\lim\limits_{\varepsilon \rightarrow 0}\, \inf\limits_{\Psi \in
X_\eps} \frac{\E\left( \Psi \right) }{\left| \HH \left( \Psi \right) \right| }
\end{equation}
where
\begin{equation} \label{icsep}
X_\eps=\left\{ \Psi \in X\ :\ \left\| \Psi \right\| _{L^{\infty }(\R^{n},V)} <\varepsilon \right\}
\end{equation}
and assuming that $m_{0}>0$, we introduce the \textit{hylomorphy ratio} $\Lambda \left( \Psi \right)$ of a state $\Psi \in X$, defined as
\begin{equation} \label{lambda}
\Lambda \left( \Psi \right) :=\frac{\E\left( \Psi \right) }{m_{0}\left| \HH\left( \Psi \right) \right| }
\end{equation}
The hylomorphy ratio of a state turns out to be a dimensionless invariant of the motion and an important quantity for the characterisation of solitons.

\begin{definition} \label{hys}
A soliton $\Psi _{0}$ is called \emph{hylomorphic} if
\begin{equation} \label{yc}
\Lambda \left( \Psi _{0}\right) <1
\end{equation}
\end{definition}

In the following we will refer to (\ref{yc}) as the \textit{hylomorphy condition}. In the study of hylomorphic solitons, the hylomorphy ratio plays a very important role and it is related to a density function which we call \textit{binding energy density}.

We assume  $\E$ and $\HH$ to be local quantities, namely, given $\Psi \in X$ there exist density functions $\rho_{\E,\Psi}\left( x\right)$ and $\rho_{\HH,\Psi }\left(
x\right) \in L^{1}(\R^{n})$ such that
\begin{eqnarray}
\E \left( \Psi \right) &=&\int \rho_{\mathcal{E,}\Psi }\left(x\right)\ dx \label{density-en}\\
\HH \left( \Psi \right) &=&\int \rho_{\mathcal{H,}\Psi }\left(x\right)\ dx \label{density-cha}
\end{eqnarray}
Then we introduce the \emph{binding energy density} defined as
\begin{equation} \label{bond-energy}
\rho_B(x)=\rho_{B,\Psi}(x) := \left[m_{0} \left| \rho_{\HH,\Psi } (x) \right| -\rho_{\E,\Psi } (x)\right]^{+}
\end{equation}
The support of $\rho_{B,\Psi}(x)$
\begin{equation} \label{support}
\Sigma \left( \Psi \right) =\overline{\left\{ x\ :\ \rho_{B,\Psi}(x)\neq 0\right\} }.
\end{equation}
is called the \emph{condensed matter region} since in these points
the binding forces prevail, and the quantity
\begin{equation}\label{conde-matter}
\int_{\Sigma \left( \Psi \right)}\rho_{\E,\Psi } (x) dx
\end{equation}
will be called the \emph{condensed matter}. The next proposition justifies the choice of the name ``hylomorphic'', since it shows that hylomorphic solitons ``contain condensed matter''.

\begin{prop} \label{supp_hylo}
If $\Lambda\left( \Psi _{0}\right) <1 $ then for all $t\in \R$ the support $\Sigma(U_{t}\Psi_{0})$ of $\rho_{B,U_{t}\Psi }(x)$ is not empty.
\end{prop}

\noindent \textbf{Proof.} It follows by the simple relations
$$
\int \rho_{B,U_{t}\Psi }(x) = \int \left[ m_{0} \left| \rho_{\HH,U_{t}\Psi } (x) \right| -\rho_{\E,U_{t}\Psi } (x) \right] ^{+} \ge
$$
$$
\ge m_{0}| \HH\left( U_{t}\Psi \right)| -\E\left( U_{t}\Psi \right)=  m_{0}|\HH \left( \Psi \right)| -\E\left( \Psi \right) =
$$
$$
=m_{0}|\HH \left( \Psi \right)| \left( 1-\Lambda \left( \Psi \right) \right) >0
$$
\qed

\bigskip

If $\Psi (x)$ is a finite energy field usually it disperses as time goes on, namely
$$
\lim\limits_{t \to \infty } \left\| U_{t}\Psi (x)\right\|_{L^{\infty }(\R^{n},V)} = 0
$$
However, if $\Lambda \left( \Psi \right) <1$ this is not the case.

\begin{prop} \label{not-dispersion}
If $\Lambda \left( \Psi \right) <1$ then
$$
\liminf\limits_{t \to \infty } \left\| U_{t}\Psi \right\|_{L^{\infty }(\R^{n},V)}=\delta >0
$$
\end{prop}

\noindent \textbf{Proof.} Let $\Lambda \left( \Psi \right) =1-a$ for a given $a>0$. We argue indirectly and assume that, for every $\varepsilon >0, $ there exists $\bar{t}$ such that
$$
\left\| U_{\bar{t}}\Psi \right\| _{L^{\infty }(\R^{n},V)}<\varepsilon
$$
namely, $U_{\bar{t}}\Psi \in X_{\varepsilon }$ where $X_{\varepsilon }$ is
defined as in (\ref{icsep}). Then, by (\ref{brutta}), if $\varepsilon $ is
sufficiently small
$$
\frac{\E \left( U_{\bar{t}}\Psi \right) }{\left| \HH\left( U_{\bar{t}}\Psi \right) \right| }\geq m_{0}-\frac{am_{0}}{2}
$$
and hence
$$
\Lambda \left( U_{\bar{t}}\Psi \right) \geq 1-\frac{a}{2}>1-a
$$
Since $\Lambda \left( U_{\bar{t}}\Psi \right) =\Lambda \left( \Psi \right)$
we get a contradiction. \qed

\bigskip

Thus if $\Lambda \left( \Psi \right) <1$, by the above propositions the
field $\Psi$ and the condensed matter  will not disperse, but will
form bumps of matter which eventually might lead to the formation of one or more hylomorphic solitons.

\section{The Nonlinear Klein-Gordon equation} \label{nonlinear-kg}

The nonlinear Klein-Gordon equation (NKG) is given by
\begin{equation} \label{KG}
\square \psi +W^{\prime }(\psi )=0  \tag{NKG}
\end{equation}
where $\square =\partial_{t}^{2}-\Delta$, $\psi(t,x): \R \times \R^{n} \to \mathbb{C}$ and $W:\mathbb{C}\rightarrow \mathbb{R}$ satisfies
$$
W(\psi)=F(|\psi|)
$$
for some smooth function $F:\mathbb{R}^{+}\rightarrow \mathbb{R}$. Also we have used the notation
\begin{equation*}
W^{\prime }(\psi )=F^{\prime }(\left| \psi \right| )\frac{\psi }{\left| \psi \right| }.
\end{equation*}
>From now on, we always assume that
\begin{equation*}
W(0)=W^{\prime}(0)=0
\end{equation*}
Under very mild assumptions on $W$, the (NKG) admits spherical solitons. In this section, we recall some general facts about (NKG) and study the hylomorphic properties of solitons. In the next section we study the problem of existence and stability of solitons, and give a classification of the nonlinear terms in (NKG) according to the properties of the solutions it admits.

We recall the pioneering paper of Rosen \cite{rosen68} and \cite{Coleman78}, \cite{strauss}, \cite{Beres-Lions}. In Physics, the spherically symmetric solitary waves have been called $Q$-balls by Coleman in \cite{Coleman86} and this is the name used in the Physics literature.

\subsection{General features of NKG} \label{gen-feat-nkg}

In this case, the fields are functions $\Psi = (\psi, \psi_{t})$ with values in $V= \C^{2}$ and equation (\ref{KG}) is the Euler-Lagrange equation of the action functional
\begin{equation} \label{az}
\mathcal{S}=\int \mathcal{L}\left( \psi_{t} ,\nabla \psi ,\psi \right) dxdt
\end{equation}
where the Lagrangian density is given by
\begin{equation} \label{lag}
\mathcal{L}=\frac{1}{2}\left| \psi_{t} \right| ^{2}-\frac{1}{2}
|\nabla \psi |^{2}-W(\psi ).
\end{equation}
Moreover the nonlinear Klein-Gordon equation admits a formulation
as an infinite dimensional Hamiltonian dynamical system. The phase
space is given by $X=H^{1}(\R^{n},\C) \times L^{2}(\R^{n},\C)$ and
the state of the system at the time $t$ is then defined by a point
$\Psi =(\psi ,\psi_{t})\in X$, where $\psi(t,x)$ is an admissible
function for the functional (\ref{az}) and $\psi _{t}:=\partial
_{t}\psi$. Equation (NKG) written as a first order system on $X$
takes the form
$$
\frac{\partial {\Psi }}{\partial t}=A{\Psi }-G({\Psi })
$$
where
\begin{equation*}
A=\left(
\begin{array}{cc}
0 & 1 \\[0.2cm]
\triangle & 0
\end{array}
\right) \qquad \qquad G({\Psi })=\left(
\begin{array}{c}
0 \\[0.2cm]
W^{\prime }(\psi )
\end{array}
\right) .
\end{equation*}
Hence, looking at $\psi_t$ as the variable conjugate to $\psi$, the dynamics on $X$ is an infinite dimensional Hamiltonian system with the Hamiltonian function given by
$$
H(\Psi)=\frac{1}{2}\left| \psi _{t}\right| ^{2}+\frac{1}{2}|\nabla \psi|^{2}+W(\psi ).
$$
The nonlinear Klein-Gordon equation is the simplest equation invariant for the Poincar\'{e} group and the action of $S^{1}$ on $X$ given by
\begin{equation} \label{july}
S^{1}\times X \ni (\theta,\Psi) \mapsto T_{\theta }\left(
\Psi\right) =\left( e^{i\theta }\psi , e^{i\theta }\psi_{t}
\right) \in X
\end{equation}
By Noether theorem the existence of conservation laws implies the existence of integrals of motion. In particular, the time invariance of the Lagrangian implies the preservation of the energy $\E$ given by
$$
\E(\psi,\psi_{t})=\int \left[ \frac{1}{2}\left| \psi_{t} \right| ^{2}+
\frac{1}{2}\left| \nabla \psi \right| ^{2}+W(\psi )\right] dx
$$
whereas the invariance with respect to the $S^{1}$-action (\ref{july}) implies the preservation of the hylomorphic charge $\HH$ given by
$$
\HH(\psi,\psi_{t}) = \int \mbox{Im} \left( \psi_{t}\, \overline{\psi }\right) \;dx
$$

In the analysis of the behaviour of solutions of (NKG), it turns useful to write $\psi$ in polar form, namely
\begin{equation}  \label{polar}
\psi (t,x)=u(t,x)e^{iS(t,x)}
\end{equation}
where $u(t,x)\in \R^{+}$ and $S(t,x)\in \R/(2\pi \Z)$. Letting $u_t:=\partial_{t}u$, $k(t,x):=\nabla S(t,x)$ and $\omega(t,x):=-\partial_{t}S(t,x)$, a state $\Psi \in X$ is uniquely determined by the quadruple $(u,u_{t},\omega ,k)$. Using these variables, the action (\ref{az}) takes the form
$$
\mathcal{S}(u,u_{t},\omega ,k)=\frac{1}{2}\int\left[ u^{2}_{t} -
\left| \nabla u\right| ^{2}+\left( \omega ^{2}-|k|^{2}\right) u^{2} \right]\, dx\, dt -
\int W(u)\, dx\, dt
$$
and equation (\ref{KG}) becomes
\begin{equation} \label{KG1}
\square u+\left( |k|^{2}-\omega ^{2}\right) u+W^{\prime }(u)=0
\end{equation}
\begin{equation}  \label{KG2}
\partial _{t}\left( \omega u^{2}\right) +\nabla \cdot \left( k u^{2}\right) =0
\end{equation}
In this form, energy and hylomorphic charge become
\begin{equation}  \label{penergy}
\E(u,u_{t},\omega ,k)=\int \left[ \frac{1}{2} u_{t}^{2}+\frac{1}{2} \left| \nabla u\right| ^{2}+\frac{1}{2}\left( \omega ^{2}+|k|^{2} \right) u^{2}+W(u)\right] dx
\end{equation}
\begin{equation} \label{cha}
\HH(u,u_{t},\omega,k)= - \int \omega \,u^{2}dx
\end{equation}

We now describe a possible interpretation of the hylomorphic
quantities such as the hylomorphic charge, the hylomorphy ratio
etc. The equation (NKG) describes how the density of hylomorphic
charge moves with the dynamics. Indeed, looking at (NKG) in polar
form, one notices that (\ref{KG2}) is the continuity equation for
the density $\rho_{\HH,\psi}$ of the hylomorphic charge $\HH$ (see
(\ref{density-cha}) and (\ref{cha})). When $\rho_{\HH,\psi}$ is
negative, it can be interpreted as the density of ``antiparticles''.
>From this point of view, equation (\ref{KG2})  describes the
conservative transport of the hylomorphic charge by a particle
field and the quantity $m_{0}$ defined in (\ref{brutta}) can be
considered as the rest energy of each particle.

In this interpretation, $\HH\left( \Psi \right)$ represents the
total number of particles counted algebraically (particles minus
antiparticles), $\frac{\E\left( \Psi \right) }{\left|
\HH\left(\Psi \right) \right| }$ represents the average energy of
each particle, and the hylomorphy ratio $\Lambda \left( \Psi
\right)$ defined in (\ref{lambda}) is a dimensionless quantity
which normalise the average energy for particle with respect to
the rest energy $m_{0}$. If $\Lambda \left( \Psi \right) >1$, the
average energy of each particle is larger than the rest energy; if
$\Lambda \left( \Psi \right) <1$, the opposite occurs and this
fact means that particles interact with one another by an
attractive force. In this interpretation, the volume occupied by
the condensed matter defined in (\ref{support}) consists of
particles tied together. Moreover the quantity $\rho_{B,\Psi}$
defined in (\ref{bond-energy}), has the dimension of the energy
and represents the binding energy of the particles. In fact,
consider for example the case of identical particles. When the
particles are free and at rest, their energy is given by the
number of particles times the energy of a particle, namely
$m_{0}|\HH\left( \Psi \right)|$, whereas the energy of the
configuration $\Psi $ is given by $\E\left( {\Psi }\right)$. Thus
the energy necessary to ``free'' the particles is $m_{0}|\HH|-\E$.
Moreover $\int \rho_{B,\Psi} (x)dx$ represents the portion of the
binding energy which is localised in the support $\Sigma(\Psi)$.

Next we will compute the hylomorphy ratio $\Lambda$ for (NKG).

\begin{theorem} \label{hr}
For (NKG) with $W$ of class $C^{2}$  the hylomorphy ratio defined in (\ref{lambda}) takes the form
$$
\Lambda \left( \Psi \right) =\frac{\E\left( \Psi \right) }{\left| \HH\left( \Psi \right) \right| \sqrt{W^{\prime \prime }(0)}}
$$
\end{theorem}

\noindent{\bf Proof.} Let us write
\begin{equation}  \label{wn}
W(s)=\frac{1}{2}m^{2}s^{2}+N(s) \qquad s\in \R^{+}
\end{equation}
where $m^{2}=W^{\prime \prime }(0)$. We have to prove that the term $m_{0}$ defined in (\ref{brutta}) satisfies $m_{0} = m$.

Using the polar form (\ref{polar}) for $\Psi$ and (\ref{penergy}), (\ref{cha}) and (\ref{wn}), we get
\begin{eqnarray*}
\frac{\E\left( \Psi \right) }{\left| \HH\left( \Psi \right) \right| } &=&\frac{\int \left[ \frac{1}{2} u_{t}^{2}+ \frac{1}{2}\left| \nabla u\right| ^{2}+\frac{1}{2}\left( \omega ^{2}+|k|^{2}\right) u^{2}+W(u)\right] dx}{\left| \int \omega \,u^{2}dx\right| } \ge \\[0.2cm]
&\geq &\frac{\int \left[ \frac{1}{2}\omega ^{2}u^{2}+\frac{1}{2} m^{2}u^{2}+N(u)\right] dx}{\int \left| \omega \right| \,u^{2}dx}
\end{eqnarray*}
Since by the classical Cauchy-Schwarz inequality
\begin{eqnarray*}
\int \left| \omega \right| \,u^{2}dx &\leq &\left( \int \omega
^{2}\,u^{2}dx\right) ^{1/2}\cdot \left( \int \,u^{2}dx\right) ^{1/2} = \\[0.2cm]
&=&\frac{1}{m}\left( \int \omega ^{2}\,u^{2}dx\right) ^{1/2}\cdot \left(
\int \,m^{2}u^{2}dx\right) ^{1/2} \le \\[0.2cm]
&\leq &\frac{1}{2m}\left[ \int \omega ^{2}u^{2}dx+\int \,m^{2}u^{2}dx\right] =
\\[0.2cm]
&=&\frac{1}{2m}\int \left( \omega ^{2}+m^{2}\right) u^{2}dx
\end{eqnarray*}
we have that
$$
\frac{\E\left( \Psi \right) }{\left| \HH\left( \Psi \right)
\right| }\geq \frac{\int \left[ \frac{1}{2}\omega ^{2}u^{2}+\frac{1}{2}m^{2}u^{2}+N(u)\right] dx}{\frac{1}{2m}\int \left( \omega ^{2}+m^{2}\right) u^{2}dx}=m+\frac{\int N(u)dx}{\frac{1}{2m}\int \left( \omega^{2}+m^{2}\right) u^{2}dx}
$$
Then since $N(u)=o(u^{2})$ for $u\rightarrow 0$, it follows that
$$
m_{0} = \lim\limits_{\eps \to 0} \inf\limits_{\Psi \in X_{\eps}}
\frac{\E\left( \Psi \right) }{\left| \HH\left( \Psi \right) \right| } \geq m
$$

In order to prove the opposite inequality, consider the family of states $\Psi _{\delta ,R}=(\delta u_{R},-i m\delta u_{R})$, where
$$
u_{R}(x)=\left\{
\begin{array}{ll}
1 & \mbox{if }|x|\le R \\[0.1cm]
0 & \mbox{if }|x|\ge R+1 \\[0.1cm]
1+R-|x| & \mbox{if }R\le |x| \le R+1
\end{array} \right.
$$
Then
\begin{eqnarray*}
\inf\limits_{\Psi \in X_{\eps}} \frac{\E\left( \Psi \right) }{\left| \HH\left( \Psi \right) \right| } &\leq &\frac{\E\left( \Psi_{\eps,R}\right) }{\left| \HH\left( \Psi_{\eps,R}\right) \right| }=\frac{\varepsilon ^{2}\int \left[ \frac{1}{2}\left| \nabla u_{R}\right| ^{2}+\frac{m^2}{2}u_{R}^{2}+\frac{1}{\eps^{2}}W(\eps u_{R})\right] dx}{\varepsilon ^{2} m\int \,u_{R}^{2}dx} = \\[0.2cm]
&=&\frac{\int \left[ \frac{1}{2}\left| \nabla u_{R}\right| ^{2}+m^2u_{R}^{2}+ \frac{1}{\eps^{2}}N(\eps u_{R})\right] dx}{m \int \,u_{R}^{2}dx} \le \\[0.2cm]
&\le & m+\frac{1}{2 m}\frac{\int \left| \nabla u_{R}\right| ^{2}dx}{\int \,u_{R}^{2}dx}+\frac{\int N(\eps u_{R})dx}{\eps^{2} m \int \,u_{R}^{2}dx} = \\[0.2cm]
&=&m+O\left( \frac{1}{R}\right) +o\left( \eps \right).
\end{eqnarray*}
where we have again used $N(u)=o(u^{2})$ for $u\rightarrow 0$. \qed

\section{Q-balls}

This section is devoted to the existence of $Q$-balls and to the study of
their structure.

\subsection{Solitary waves for NKG} \label{solwavNKG}

The easiest way to produce solitary waves of (\ref{KG}) consists in solving
the static equation
\begin{equation} \label{KGS}
-\Delta u+W^{\prime }(u)=0
\end{equation}
and making a change of the frame of reference to give a velocity to the wave, setting for example
\begin{equation}  \label{pina}
\psi_{v}(t,x)=\psi_{v}(t,x_{1},\dots,x_{n})=u\left( \frac{x_{1}-vt}{\sqrt{1-v^{2}}},\dots,x_{n}\right)
\end{equation}
which is a solution of (\ref{KG}) which represents a bump travelling in the $x_{1}$-direction with speed $v$ (we consider a normalisation of units of measure so that the speed of light $c$ is equal to 1).

Unfortunately Derrick, in the very well known paper  \cite{D64}, has proved
that the request $W(u)\geq 0$ (which is necessary if we want the energy to be a non-negative invariant) implies that equation (\ref{KGS}) has only the trivial solution. His proof
is based on the Derrick-Pohozaev identity, which for (\ref{KGS}) is given by
$$
\left( \frac{1}{n}-\frac 12 \right)\int \left| \nabla u\right| ^{2}dx\, = \int W(u)\, dx
$$
and holds for any finite energy solution $u$ of equation (\ref{KGS}) (for details
see also \cite{sammomme}). Clearly the above equality for $n\ge 3$ and $W(u)\geq 0$ implies that $u\equiv 0$.

However, we can try to prove the existence of solitons for (\ref{KG}) exploiting the possible existence of \textit{standing waves}, namely finite energy solution of (\ref{KG}) of the form
\begin{equation} \label{sw}
\psi _{0}(t,x)=u_{0} (x)e^{-i\omega _{0}t}, \qquad u_{0} \geq 0.
\end{equation}
Substituting (\ref{sw}) in (\ref{KG}), one finds that $u_{0}(x)$ is a solution of
\begin{equation}  \label{static}
-\Delta u+W^{\prime }(u)=\omega _{0}^{2}u
\end{equation}
and the existence of non-trivial solutions of (\ref{static}) is not prevented by the Derrick-Pohozaev identity, which in this case reads
\begin{equation} \label{ventuno1}
\left( \frac{1}{n}-\frac 12 \right)\int \left| \nabla u\right| ^{2}dx\, = \int \left(W(u) - \frac 1 2 \omega_{0}^{2} u^{2}\right)\, dx
\end{equation}

Since the Lagrangian (\ref{lag}) is invariant for the Lorentz group, we can
obtain other solutions $\psi _{\mathbf{v}}(t,x)$ just making a Lorentz
transformation on it. Namely for $n=3$, if we take the velocity $\mathbf{v}=(v,0,0),$ $\left| v\right| <1$, and set
\begin{equation*}
t^{\prime }=\gamma \left( t-vx_{1}\right) ,\text{ }x_{1}^{\prime }=\gamma
\left( x_{1}-vt\right) ,\text{ }x_{2}^{\prime }=x_{2},\text{ }x_{3}^{\prime
}=x_{3}
\end{equation*}
with
\begin{equation*}
\gamma =\frac{1}{\sqrt{1-v^{2}}},
\end{equation*}
it turns out that
\begin{equation*}
\psi _{\mathbf{v}}(t,x):=\psi (t^{\prime },x^{\prime })
\end{equation*}
is a solution of (\ref{KG}).

In particular, given a standing wave $\psi_{0} (t,x)=u_{0}(x)e^{-i\omega _{0}t},$
the function $\psi_{\mathbf{v}}(t,x)=\psi_{0} (t^{\prime },x^{\prime })$ is a
solitary wave which travels with velocity $\mathbf{v}$. Thus, if for example $\mathbf{v}= (v,0,0)$ and $u(x)=u(x_{1},x_{2},x_{3})$ is any solution of equation (\ref{static}), then
\begin{equation} \label{solitone}
\psi _{\mathbf{v}}(t,x_{1},x_{2},x_{3})=u \left( \gamma \left(x_{1}-vt\right) ,x_{2},x_{3}\right) e^{i(\mathbf{k\cdot x}-\omega t)},
\end{equation}
is a solution of (\ref{KG}) provided that
$$
\omega =\gamma \omega_{0} \quad \text{and}\quad \mathbf{k}=\gamma \omega_{0}\mathbf{v}
$$
Notice that (\ref{pina}) is a particular case of (\ref{solitone}) when $\omega _{0}=0$.

Finally we remark that if $u$ is a solution of equation (\ref{static}), then by a well known result on elliptic equations \cite{GiNiNi} it has necessarily spherical symmetry. Coleman called the spherically symmetric solitary waves of (NKG) \emph{Q-balls} (\cite{Coleman86}) and this is the name generally used in the Physics literature.

In order to prove the existence of solitons for (NKG), we are now led to study the existence of orbitally stable standing waves. Hence in particular we study existence results of couples $(u_{0},\omega_{0})$ which satisfy equation (\ref{static}) under some general assumptions on the function $W$.

\subsection{Existence results for Q-balls} \label{qball}

From now on we make the following assumptions on $W$:

\begin{itemize}
\item  (W-i) \textbf{(Positivity)} $W(s)\geq 0$;

\item  (W-ii) \textbf{(Normalisation)} $W(0)=W^{\prime }(0)=0;\;W^{\prime
\prime }(0)=1$;

\item  (W-iii) \textbf{(Hylomorphy)} $\lambda _{0}:=\inf\limits_{s\in \R^{+}} \frac{W(s)}{\frac{1}{2} s^{2}}<1$.
\end{itemize}

Let us make some remarks on these assumptions.

\bigskip

(W-i) implies that the energy is positive for any state. Some aspects of the theory remain true weakening this assumption. In particular there are results in the case $n\ge 3$ also for
\begin{equation} \label{oc}
W(s)=\frac{1}{2} s^{2}-\frac{1}{p} s^{p}, \qquad 2<p<\frac{2n}{n-2},
\end{equation}
since they are easier to prove with variational methods. However assumption (W-i) is required by most of the physical models and it simplifies some theorems.

\bigskip

(W-ii) is a normalisation condition. In order to have solitary waves it is
necessary that $W^{\prime \prime }(0)\geq 0.$ There are results also
in the null-mass case $W^{\prime \prime }(0)=0$ (see e.g. \cite{Beres-Lions}
and \cite{BBR07}). However the most interesting situations occur when $W^{\prime \prime }(0)>0.$ In the latter case, re-scaling the independent variables $(t,x)$, we may assume that $W^{\prime \prime }(0)=1$.

\bigskip

(W-iii) is the crucial assumption. As we will see this assumption is necessary to have states $\Psi $ with $\Lambda \left( \Psi \right) <1$, hence hylomorphic solitons according to Definition \ref{hys}.

\bigskip

Under the previous assumptions, we can write $W$ as
\begin{equation} \label{NN}
W(s)=\frac{1}{2}s^{2}+N(s), \qquad s\in \R^{+}
\end{equation}
and the hylomorphy condition can be stated by saying that there exists $s_{0} \in \R^{+}$ such that $N(s_{0})<0$. Actually, by our interpretation of (NKG) (see Section \ref{gen-feat-nkg}), $N$ is the nonlinear term which, when negative, produces an ``attractive interaction'' among the particles.

In \cite{BBBM} we prove that

\begin{theorem}[\cite{BBBM}] \label{main-theorem}
If (W-i), (W-ii) and (W-iii) hold then (NKG) in $\R^{n}$, with $n\ge 2$, admits hylomorphic
solitons of the form (\ref{static}).
\end{theorem}

This result is recent in the form given here. In fact only recently it has been proved the orbital stability of the standing waves (\ref{sw}) with respect to the standard topology of $X=H^{1}(\R^{n},\C)\times L^{2}(\R^{n},\C)$ and for all the $W$ which satisfy (W-ii), (W-iii) and a condition weaker than (W-i). Nevertheless Theorem \ref{main-theorem} has a very long history starting with the pioneering paper of Rosen \cite{rosen68}. Coleman \cite{Coleman78} and Strauss \cite{strauss} gave the first rigourous proofs of existence of
solutions of the type (\ref{sw}) for some particular $W$, and later necessary and sufficient existence conditions have been found by Berestycki and Lions \cite{Beres-Lions}.

The first orbital stability results for (NKG) are due to Shatah, Grillakis and Strauss \cite{shatah}, \cite{gss87}. Namely, they consider the real function
$$
\omega \mapsto d(\omega):= \E (u_{\omega},\omega) + \omega |\HH(u_{\omega},\omega)|
$$
where $u_{\omega}$ is the ground state solution of (\ref{static}) for a
fixed $\omega$, and $\E(u_{\omega},\omega)$ and $\HH(u_{\omega},\omega)$ are the energy and the hylomorphic charge of the standing wave $\psi_{\omega}(t,x)=u_{\omega}(x)e^{-i\omega t}$. They prove that a necessary and sufficient condition for the orbital stability of the standing wave $\psi_{\omega_{0}}$ is the convexity of the map
$d(\omega)$ in $\omega_{0}$. It is difficult to verify theoretically this condition for a general $W$ since $d(\omega)$ cannot be computed explicitly. In some particular cases the properties of $d(\omega)$ can be investigated, see for instance \cite{ss85} where the authors give exact ranges of the frequency $\omega$ for which they obtain stability and instability of the respective standing waves for the nonlinear wave equation with $W(\psi)$ as in (\ref{oc}). Moreover the computation of  $d(\omega)$ in a neighbourhood of a chosen $\omega_0$ can be extremely difficult due to the fact that the knowledge of $u_{\omega}$ is required.

The proof of Theorem \ref{main-theorem} follows from some interesting preliminary results, which give the idea of the proof. We now sketch the main steps, following the general method exposed in Section \ref{def-hylo-sol}, referring to \cite{BBBM} for details.

The first step is to define the following set
$$
Y := H^{1}(\R^{n},\R^{+}) \times \R
$$
which corresponds to a subset of $X$ through the identification
$$
Y \ni (u(x),\omega) \longmapsto ( u(x),-i\omega u(x)) \in X =
H^{1}(\R^{n},\C) \times L^{2}(\R^{n},\C)
$$
Notice that the standing waves are contained in $Y$.

We introduce the functionals energy and charge on the space $Y$ which with abuse of notation we still denote by $\E$ and $\HH$. They take the form
\begin{equation} \label{energy_Q_ball}
\E (u,\omega) =\int \left[ \frac{1}{2}\left| \nabla u \right|^{2} + \frac{1}{2}\omega^{2} u^{2} + W(u)\right] dx
\end{equation}
\begin{equation} \label{charge_Q_ball}
\HH (u,\omega) =-\omega \int u^{2}\, dx
\end{equation}
Without loss of generality we restrict ourselves to the case $\HH<0$ and $\omega>0$. The relationship between the (NKG) equation and the energy $\E$ on the space $Y$ is given by the following proposition
\begin{prop}[\cite{BBBM}] \label{bello}
The function $\psi (t,x)=u(x)e^{-i\omega t}$ is a standing wave for (NKG) if and only if $(u,\omega)$ is a critical point (with $u\neq 0$) of the functional $\E (u,\omega)$ constrained to the manifold
$$
\mathcal{M}_{\sigma }= \set{ \left(u,\omega \right) \in Y\ :\ |\HH\left( u,\omega \right)| =\sigma }, \quad \sigma \neq 0
$$
\end{prop}

\noindent \textbf{Proof.}
A point $(u,\omega)\in Y$ is critical for $\E(u,\omega)$ constrained to the manifold $
M_\sigma$ if and only if there exists $\lambda$ such that
\begin{equation}  \label{Lagrvinc}
\begin{array}{c}
\frac{\partial \E(u,\omega)}{\partial u}= \lambda\ \frac{\partial \HH(u,\omega)
}{\partial u} \\[0.3cm]
\frac{\partial \E(u,\omega)}{\partial \omega}= \lambda\ \frac{\partial
\HH(u,\omega)}{\partial \omega}
\end{array}
\end{equation}
By definition of energy and charge (equations (\ref{energy_Q_ball}) and (\ref
{charge_Q_ball})), equation (\ref{Lagrvinc}) can be written as
\begin{equation*}
\begin{array}{rcl}
-\Delta u+ W^{\prime}(u)+\omega^2u & = & -2\lambda\, \omega\, u \\[0.3cm]
\int\omega\, u^2dx & = & -\lambda\int u^2dx.
\end{array}
\end{equation*}
Since $\sigma \not= 0$, from the second equation $\lambda=-\omega$, and the first
becomes equation (\ref{static}). \qed

\vskip 0.5cm

This result gives a simple criterion to obtain standing waves for (NKG). Moreover we are interested in standing waves which are orbitally stable (see Definition \ref{dos}), and to this aim we look for points of minimum of $\E$ on $\mathcal{M}_{\sigma}$.  To this aim we use the notion of hylomorphy ratio, which on $Y$ takes the form
\begin{equation} \label{hyl}
\Lambda (u,\omega) =\frac{\E (u,\omega)}{|\HH (u,\omega)|}= \frac 1 2 \, \omega
+ \frac{1}{2\omega} \, \alpha(u)
\end{equation}
where
\begin{equation}  \label{alpha}
\alpha(u):= \frac{\int \left( \frac{1}{2} |\nabla u|^2 + W(u) \right)dx}{\int \frac 1 2 u^2 \, dx}.
\end{equation}
One immediately gets
\begin{equation}  \label{lambda-alpha}
\inf\limits_{\omega \in \R^+} \, \Lambda(u,\omega) = \sqrt{\alpha(u)}
\end{equation}
and if $\lambda_0$ is defined as in (W-iii), then it is proved in \cite{BBBM} that
\begin{equation}  \label{serve}
\lambda_0 = \inf\limits_{u\in H^{1}} \, \alpha(u)
\end{equation}
To prove the existence of point of minimum for the energy we use
\begin{prop}[\cite{BBBM}] \label{lemma-existence}
If there exists $(\bar u, \bar \omega) \in  Y$ such that $\Lambda(\bar u, \bar \omega) <1$, there exist points of minimum of the energy $\E(u,\omega)$ constrained to the manifold $\mathcal{M}_{\sigma }$ with $\sigma=|\HH(\bar u,\bar \omega)|$.
\end{prop}

Equations (\ref{lambda-alpha}) and (\ref{serve}) imply that assumption (W-iii) is sufficient for the existence of $(u,\omega)$ with $\Lambda(u,\omega) <1$. Hence step (S-1) of Section \ref{def-hylo-sol} is completed.

A point of minimum $(u_{0},\omega_{0}) \in Y$ corresponds to the standing wave $\psi_{0}(t,x) = u_{0}(x) e^{-i\omega_{0}t}$, with $u_{0}(x)$ a spherically symmetric function.  Steps (S-2) and (S-3) follow from the fact that for isolated points of minimum the set of minimisers $\Gamma$ consists of the set
{\small
$$
\Gamma:= \left\{\left(u_{0}(x+a)\, e^{i\theta}, -i\omega_{0}
u_{0}(x+a)\, e^{i\theta}\right) \right\}_{a\in \R^{n},\, \theta
\in \R} \subset X
$$}
which has dimension $n+1$. In the following we restrict ourselves to the case of isolated points of minimum, which is ``generic'' in the family of (NKG) equations.

To finish the proof of Theorem \ref{main-theorem} it remains to
prove step (S-4), that is the stability of the set $\Gamma$. In
\cite{BBBM} we prove that the function
$$
V(\Psi) := \left( \E(\Psi) - \min_{\mathcal{M}_{\sigma }} \E \right)^{2} + \left( \HH(\Psi) - \sigma \right)^{2}
$$
is a Lyapunov function on $X$. This implies that
\begin{theorem}[\cite{BBBM}] \label{thm-stability}
If $(u_{0},\omega_{0})$ is a point of local minimum of the functional $\E(u,\omega)$ constrained to the manifold $M_{\sigma}$ with $\sigma=|\mathcal{H}(u_{0},\omega_{0})|$, then $\psi_{0}(t,x)=u_{0}(x)e^{-i\omega_{0} t}$ is an orbitally stable standing wave.
\end{theorem}

Theorem \ref{main-theorem} is a pure existence result and gives no
information on the charge and frequency of the standing waves. In
\cite{BBBM} we show that Proposition \ref{lemma-existence} implies
that, under the assumptions (W-i), (W-ii) and (W-iii), there
exists a threshold value $\sigma_0>0$ such that for any $\sigma
\in (\sigma_0,\infty)$ there are hylomorphic solitons with
hylomorphic charge $\sigma$. The existence of hylomorphic solitons
for all hylomorphic charges can be obtained by a stronger version
of assumption (W-iii). Let us consider the condition
\begin{itemize}
\item (W-iv) \textbf{(Behaviour at $s=0$)} $N(s)<0$ for $s\in
\R^{+}$ small enough ($N(s)$ is defined in (\ref{NN})).
\end{itemize}
\begin{cor}[\cite{BBBM}] \label{123-carica-tutte}
If (W-i), (W-ii) and (W-iv) hold then for any $\sigma \not= 0$ there exists a hylomorphic soliton for (NKG) with hylomorphic charge $|\HH|=\sigma$.
\end{cor}

We finish this section by giving a remark which is useful for the numerical approach to the existence of hylomorphic solitons for (NKG). We have found solitons as points of minimum for the energy $\E(u,\omega)$ on the manifold with fixed hylomorphic charge. Hence we are studying a minimisation problem in two variables with one constraint. It is immediate that this problem can be translated into a minimisation problem in one single variable with no constraints. We will use as independent variable the functions $u\in H^{1}$.

If $\left(u_{0},\omega _{0}\right) $ is a minimiser of $\mathcal{E}\left(u,\omega \right)$ constrained to $\mathcal{M}_{\sigma}$ with $\sigma$ fixed, then it is also a minimiser of $\Lambda \left( u,\omega \right) $ constrained to $\mathcal{M}_{\sigma }$. Moreover if $\left( u,\omega \right) \in \mathcal{M}_{\sigma }$, then
\begin{equation} \label{omega_dalla_carica}
\omega =\omega_{\sigma} \left( u \right) :=\frac{\sigma}{\int u^{2}\, dx}
\end{equation}
Hence, letting
\begin{equation} \label{ganzetto}
\Lambda_{\sigma} \left(u\right) := \Lambda \left( u,\omega_{\sigma} \left( u
\right) \right) =\frac{1}{\sigma }\int \left[ \frac{1}{2}\left| \nabla
u\right| ^{2}+W(u)\right] dx +\frac{\sigma }{2\int u^{2}dx}
\end{equation}
we can state Theorem \ref{thm-stability} in the form
\begin{theorem} \label{ganzo}
For any fixed $\sigma \in \R^{+}$, if $u_{0}(x)$ is a point of minimum for the functional $\Lambda_{\sigma} \left(u\right)$ defined in (\ref{ganzetto}), then $\psi_{0}(t,x)=u_{0}(x)e^{-i\omega_{\sigma}(u)t}$ is an orbitally stable standing wave of (NKG) with hylomorphic charge $|\HH|=\sigma $ and $\omega_{\sigma}(u)$ given as in (\ref{omega_dalla_carica}).
\end{theorem}

\subsection{Numerical construction of Q-balls} \label{minimizzazione_numerica}

In this section we provide basic details of our numerical method for constructing hylomorphic solitons.
The numerical results are shown in the next section, where they illustrate a classification of (NKG) equations originally
introduced therein.

Theorem \ref{ganzo} is straightforwardly exploited for the numerical construction of Q-balls:
we fix a hylomorphic charge $\sigma \in \R^{+}$ and look for points of minimum of the functional $\Lambda_{\sigma}(u)$.
Firstly, by the classical principle of symmetric criticality (see \cite{Pa79}), we can restrict ourselves to the analysis of radially symmetric functions $u(x)=u\left(r\left(x\right)\right)$, with $r(x)=|x|$ in $\R^{n}$.
We then consider the following evolutionary problem:
\begin{equation} \label{flusso_parabolico}
\left\{
\begin{array}{rclcl}
    \partial_t u (r,t) & = &  - \sigma \, \mathrm{d} \Lambda _{\sigma} = \Delta u - W^\prime(u) + \omega^2 u & \mathrm{in} &
    [0,\tilde{r})\, \times \R^+ \\
    u(r,t) & = & 0 & \mathrm{on} & \set{r = \tilde{r}} \times \R^+
\end{array}
\right.
\end{equation}
in which $\tilde{r}$ denotes a chosen upper bound for the $r-$domain (discussed below) and
$\omega=\omega_{\sigma}(u)$ as in (\ref{omega_dalla_carica}).
The evolution of $u$ according to (\ref{flusso_parabolico}) is a gradient flow and therefore
a non-increasing trend for $\Lambda_{\sigma} (t):= \Lambda_{\sigma} (u(r,t))$ is obtained, well suited as
the sought minimisation process.

The problem (\ref{flusso_parabolico}) is then discretised by a classical line method; in particular,
2nd order and 1st order accurate finite differences have been respectively used for space and time
discretisation (see e.g. \cite{quart}).
The chosen charge $\sigma$ is directly enforced at the $n$-th time level ($n=0,1,2,\dots$), by
evaluating the frequency $\omega^n$ from the corresponding
numerical solution through the discrete counterpart of
(\ref{omega_dalla_carica}). Moreover, time-advancing is stopped when
$|\omega^{n+1}-\omega^n|/\omega^n < e_\omega$, $e_\omega$ being a predefined
threshold (a relative error on $\Lambda_\sigma$ has been considered as well).
The proposed method manages to efficiently converge by starting from several initial guesses: not only
from Gaussian (Q-ball like) profiles but also from discontinuous ones (e.g.
$\bar{u} \cdot \chi \left( D_r \right)$, with $D_r$ contained within the
chosen $r-$domain and $\bar{u}$ suitably set for obtaining the desired $\sigma$).
Finally, the domain extreme $\tilde{r}$ is chosen in such a way that it does not affect
the numerical results (an {\it a-posteriori} check might be necessary: the chosen
domain must be large enough to contain the soliton support, with some margin).

We remark that we did not try to implement a shooting method, on purpose.
Indeed, sign-changing solutions are known to be unstable and, based on the given
definition of soliton (implying stability), such solutions are of no interest
in the present study. Conversely, the proposed numerical method
guarantees to find a soliton, due to Theorem \ref{ganzo}.

\begin{figure}
    \begin{center}
    \subfigure[]
    {\includegraphics[width=6.0cm,height=4.5cm]{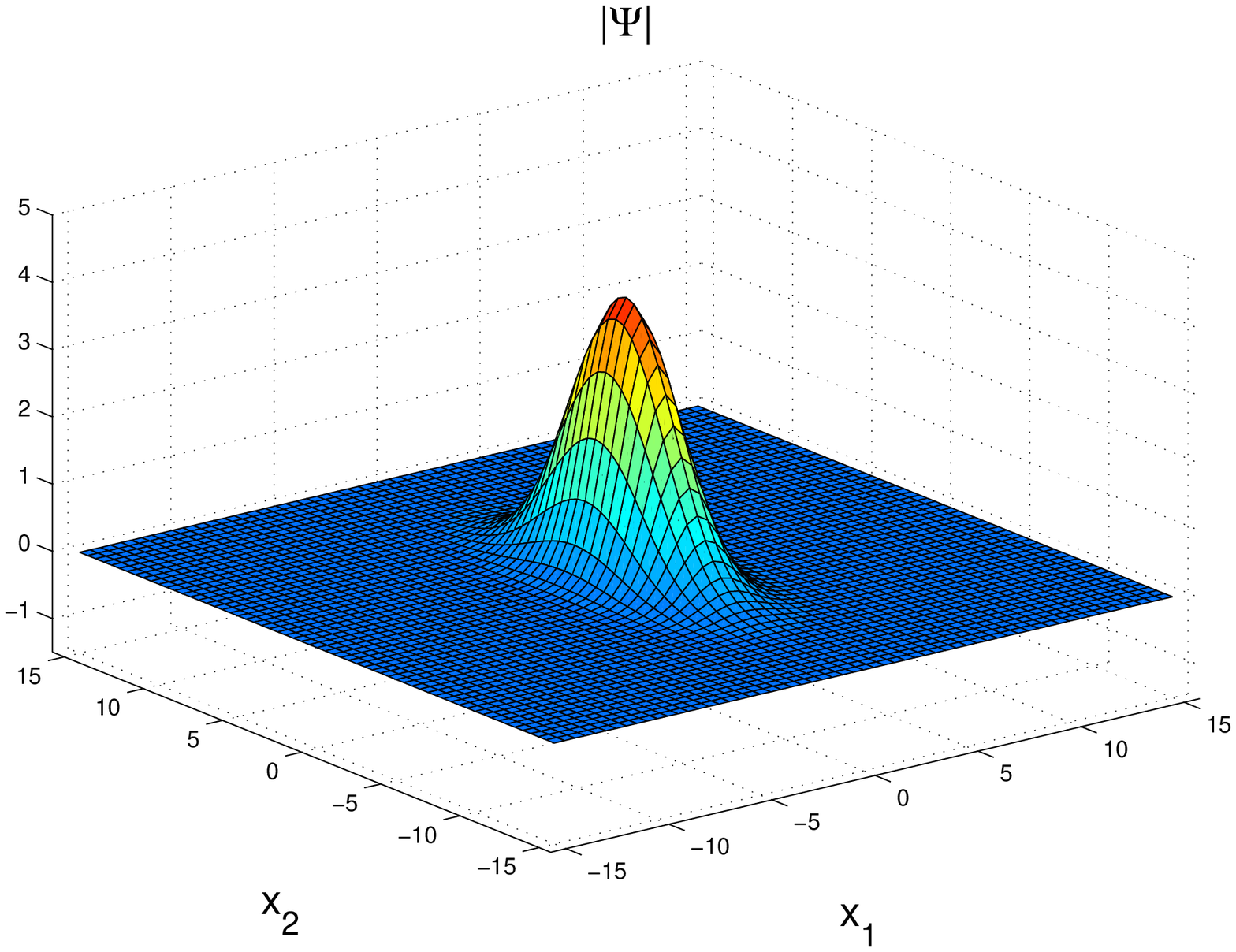}}
    \subfigure[]
    {\includegraphics[width=6.0cm,height=4.5cm]{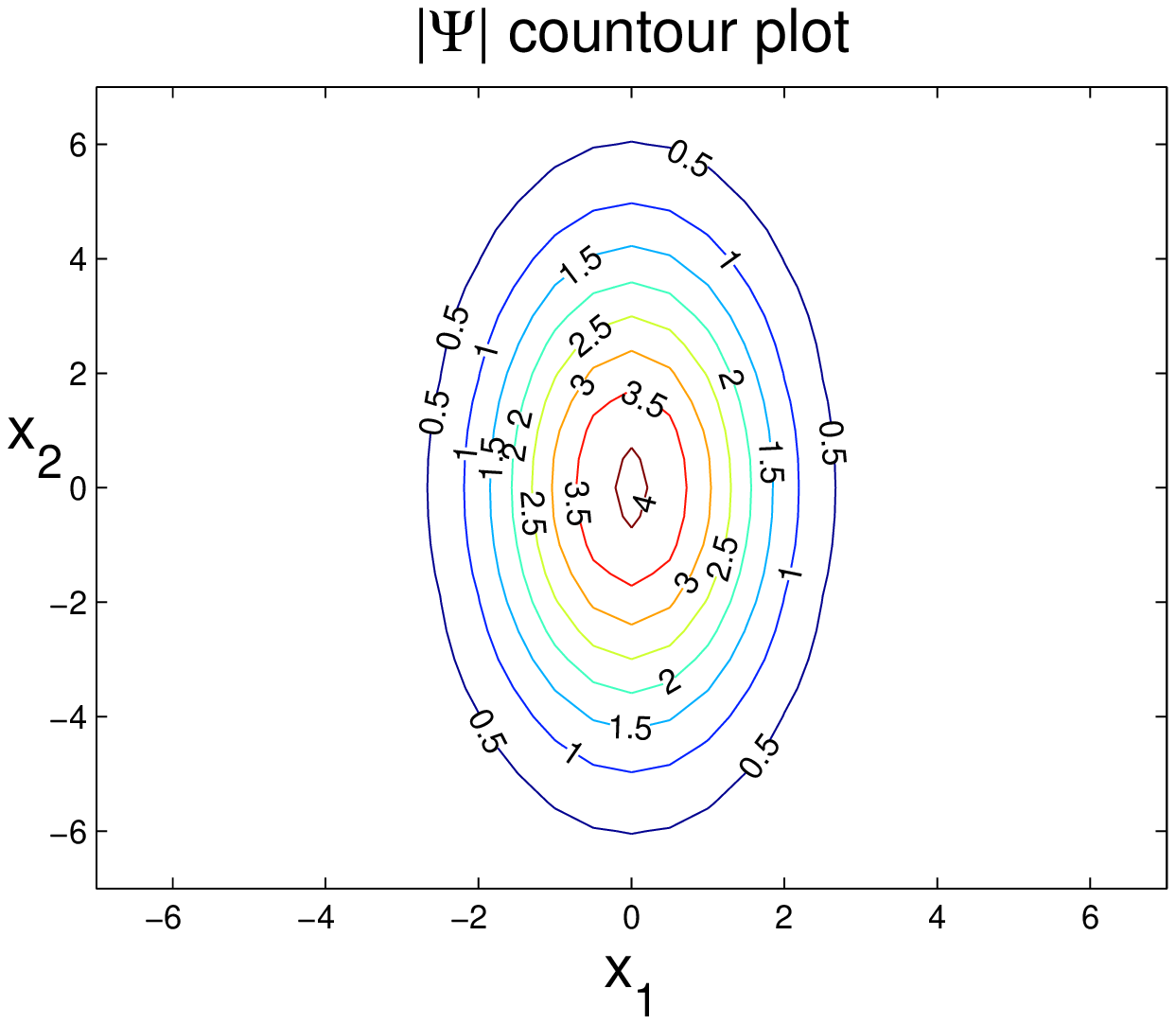}}
    \caption{(a) Surface plot of $|\Psi|$ for a 2D Q-ball moving with velocity $\mathbf{v}=(0.9,0)$
                 on the $\left(x_1,x_2\right)$ plane (recall that the units of measure are scaled
                 so that the speed of light is $c=1$). $W$ given by (\ref{ex_W_gamma}), hylomorphic
                 charge $\sigma=300$.
             (b) Corresponding contour plot (detail), highlighting
                 the Lorentz contraction along the direction of motion.}
    \label{fig_2DQball}
    \end{center}
\end{figure}

It is worth remarking that, once defined $u(r(x))$ by means of the aforementioned strategy, it is easy to build moving Q-balls through the transformation (\ref{solitone}); a
two-dimensional (2D) example is shown in Figure \ref{fig_2DQball}.

%%%%%%%%%%%%%%%%%%%%%%%%%%%%%%%%%%%%%%%%%%%%%%%%%%%%%%%%%%%%%%%%%%%%%%%%%%%%%

\subsection{Classification of the nonlinear terms} \label{classification}

In this section we introduce four classes of behaviour for the nonlinear term $W$, all classes admitting hylomorphic solitons. The classification is based on the existence and non-existence of solitons with small charge $\HH$ and with big $L^{\infty}$ norm for the (NKG) with $W$ in a fixed class. Moreover we state some general properties shared by hylomorphic solitons. Without any loss of generality we consider the case $\HH<0$ and $\omega>0$.

\begin{theorem}[Admissible frequencies]  \label{freq-adm}
Let (W-i), (W-ii) and (W-iii) hold and $\psi(x,t)=u(x)e^{-i\omega t}$ be a hylomorphic soliton for (NKG) in $\R^{n}$. Set
$$
\Gamma(u) := \frac{\int|\nabla u|^2dx }{\int u^2 dx}
$$
and
\begin{equation*}
\omega_1 := \max \set{  1-\sqrt{1-\lambda_0},\ \frac 12-\frac 12 \sqrt{1-\frac{4\Gamma(u)}{n}} }
\end{equation*}
where $\lambda_0$ is defined as in (W-iii). Then $\Gamma(u) < \frac{n}{4}$ and there exists $\eta >0$ such that the frequency $\omega$ satisfies
\begin{equation*}
\omega \in \left( \omega_1,\ \frac 12+ \frac 12 \sqrt{1-\frac{4\Gamma(u)}{n}}\right) \ \subset \ \left(\frac 1 2\, \lambda_{0}, \ 1- \eta\right)
\end{equation*}
\end{theorem}

\noindent {\bf Proof.} We first show that $\Gamma <\frac{n}{4}$ and
$\omega <\frac 12+\frac 12 \sqrt{1-\frac{4\Gamma}{n}}$.

We recall that if $\psi(t,x)$ is a hylomorphic soliton for (NKG), then its radial part $u(x)$ satisfies equation (\ref{static}). Hence using the Derrick-Pohozaev identity (\ref{ventuno1}) for $u$ it follows that the energy (\ref{energy_Q_ball}) and the hylomorphy ratio (\ref{hyl}) can be written as
$$
\E(u,\omega) = \int \frac {1}{n} |\nabla u|^2 +\omega^2 u^2dx
$$
$$
\Lambda(u, \omega) = \frac{\Gamma(u)}{n\, \omega}+\omega.
$$
By imposing that $\Lambda(u,\omega)<1$, a simple computation implies the first part of the theorem.

It remains to prove that $\omega>\omega_1$. By Proposition \ref{serve} and (\ref{hyl}) any hylomorphic soliton fulfils
$$
\Lambda(u,\omega)= \frac{\alpha(u)}{2\omega}+\frac{\omega}{2}\geq \frac{\lambda_0}{2\omega}+\frac{\omega}{2}
$$
Again by imposing that $\frac{\lambda_0}{2\omega}+\frac{\omega}{2}<1$
we obtain that $\omega>\omega_1$. \qed

\bigskip

\begin{figure}
    \begin{center}
    \subfigure[]
    {\includegraphics[width=6.0cm,height=4.5cm]{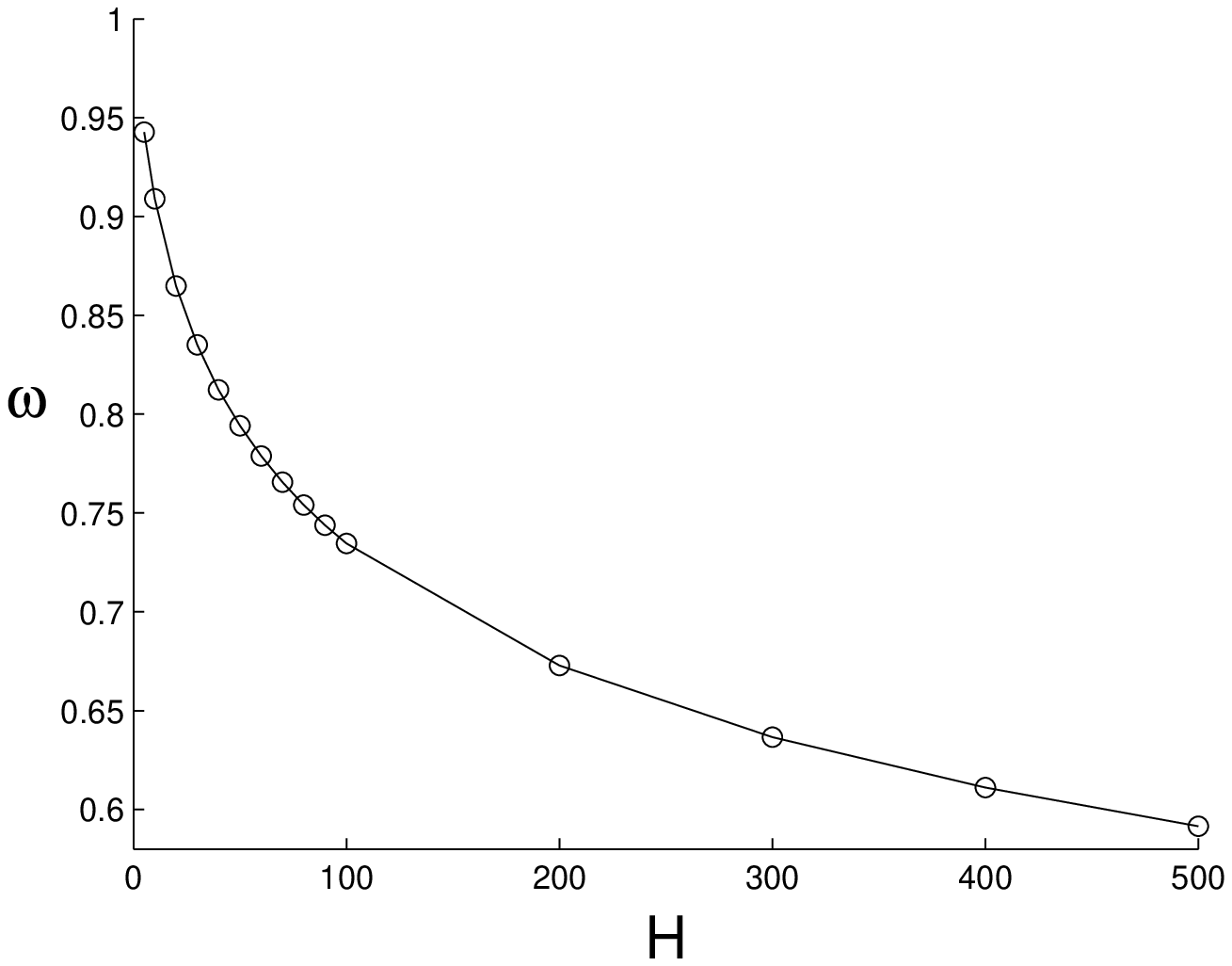}}
    \subfigure[]
    {\includegraphics[width=6.0cm,height=4.5cm]{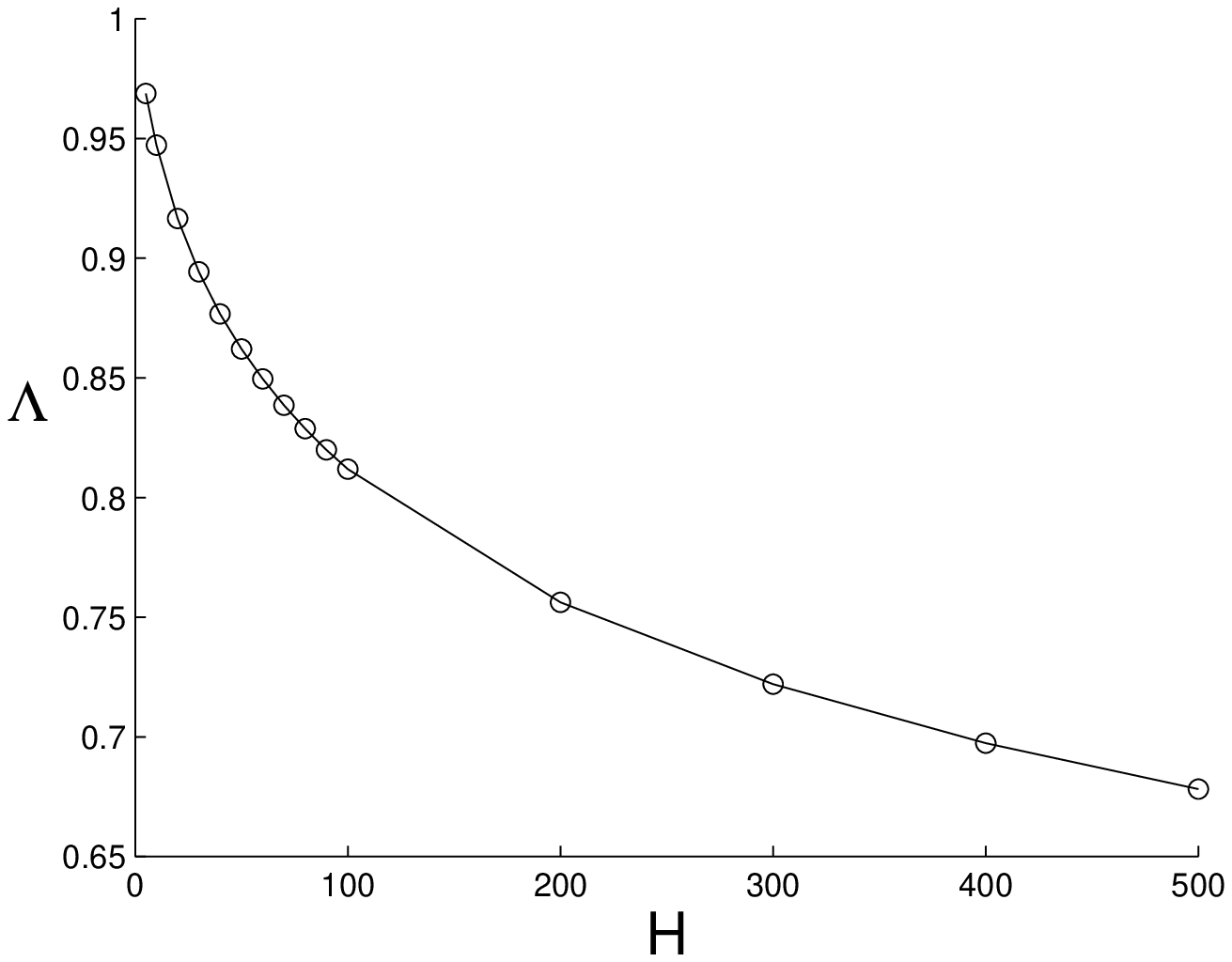}}
    \caption{Frequency (a) and hylomorphy ratio (b) of 2D Q-balls as a function of the hylomorphic charge $\HH$.
           $W$ given by (\ref{ex_W_gamma}); circles associated with the following charge values:
             $\{5,10,20,30,40,50,60,70,80,90,100,200,300,400,500\}$.}
    \label{fig_omega_e_Lambda}
    \end{center}
\end{figure}

Figure \ref{fig_omega_e_Lambda} shows a ``typical'' behaviour of the frequency and the hylomorphy ratio of a hylomorphic soliton as the hylomorphic charge varies. The frequency $\omega$ and the hylomorphy ratio $\Lambda$ are decreasing functions of the hylomorphic charge. We remark that as expected we find $\Lambda < 1$ for all hylomorphic solitons.

\begin{definition}
Let $W$ satisfy (W-i), (W-ii) and (W-iii). We say that equation (NKG) (or W) is of
\emph{type ($\alpha $)} if there exists $\alpha _{0}>0$ such that any
hylomorphic soliton $\psi (t,x)$ of the form (\ref{sw}) satisfies
$$
\left\| \psi(t,\cdot )\right\|_{L^{\infty }(\R^{n})} \ge \alpha_{0} \qquad \forall\, t \in \R
$$
\end{definition}

\begin{theorem}
Assume that $W$ satisfies (W-i), (W-ii) and (W-iii), and write it as in (\ref{NN}). If
$$
N(s)>0 \quad \forall s\in \left( 0,\alpha _{0}\right)
$$
then equation (NKG) is of type ($\alpha$).
\end{theorem}

\noindent {\bf Proof.} Let $\psi(t,x)=u(x)e^{-i\omega t}$ be a hylomorphic soliton for (NKG) and let us assume that there exists $t_{0}$ such that $\left\| \psi(t_{0},\cdot )\right\|_{L^{\infty }(\R^{n})} < \alpha_{0}$. By the assumption on $N$, this implies that $N(\psi(t_{0},x)) >0$ for all $x\in \R^{n}$.

However, by Proposition \ref{supp_hylo}, it follows that
hylomorphic solitons have non-empty support of binding energy
density at any time. That is, writing
\begin{eqnarray*}
\rho_{\mathcal{E},\psi}(t,x) &=& \frac{1}{2}\left| \nabla u\right| ^{2}+ \frac{1}{2}\left[ 1+\omega ^{2}\right] u^{2}+N(u) \\
\rho_{\HH,\psi}(t,x) &=&-\omega\, u^{2}
\end{eqnarray*}
see equations (\ref{bond-energy}) and (\ref{density-cha}), for any
$t\in \R$ there exists a positive measure set of points $x\in
\R^{n}$ such that the binding energy density fulfils
\begin{equation} \label{mu-utile}
\rho_{B,\psi} = \omega u^{2}-\frac{1}{2}\left| \nabla u\right| ^{2}-\frac{1}{2}\left[ 1+\omega
^{2}\right] u^{2}-N(u) >0
\end{equation}
see  (\ref{conde-matter}) and Theorem \ref{hr} which implies $m_{0}=1$. Moreover by equation (\ref{mu-utile}) it follows by simple computations that for all $t\in \R$ there exists a positive measure set of points $x\in \R^{n}$ such that
$$
0< \rho_{B,\psi}(t,x) \le \left[ \omega -\frac{1}{2}-\frac{1}{2}\omega ^{2}\right] u^{2}(x) -N(u(x)) \le -N(u(x))
$$
Letting $t=t_{0}$ in the previous inequality we get that $N(\psi(t_{0},x)) <0$ for some $x$, which is a contradiction. \qed

\begin{definition}
Let $W$ satisfy (W-i), (W-ii) and (W-iii). We say that equation (NKG) (or W) is of
\emph{type (non-$\alpha $)} if for any $\alpha _{0}>0$  there exists a hylomorphic soliton
$\psi (t,x)$ of the form (\ref{sw}) for which
$$
\left\| \psi(t,\cdot )\right\|_{L^{\infty }(\R^{n})} \le \alpha_{0} \qquad \forall\, t \in \R
$$
\end{definition}

\begin{theorem}
Assume that $W$ satisfies (W-i), (W-ii) and (W-iv), and that $\lambda_0$ defined in (W-iii) is positive. Then equation (NKG) is of type ($non$-$\alpha$).
\end{theorem}

\noindent {\bf Proof.} By Corollary \ref{123-carica-tutte} if (W-i), (W-ii) and (W-iv) for $W$ hold, there exist hylomorphic solitons of the form (\ref{sw}) for arbitrary small hylomorphic charges. Moreover by Theorem \ref{freq-adm}, if $\lambda_{0}>0$ the admissible frequencies of the hylomorphic solitons satisfy
$$
\omega \ge 1-\sqrt{1-\lambda_0} >0
$$
Hence we can find a sequence of hylomorphic solitons $\psi_k=u_k(x)e^{-i\omega_k t}$ for (NKG) with hylomorphic charge vanishing and satisfying
$$
\int |\nabla u_k|^2\, dx \rightarrow 0 \qquad \qquad
\int u_k^2\, dx \rightarrow 0
$$
We recall that the $(u_{k})$ are radially symmetric functions, hence there exist positive constants $C$ and $\rho$, only depending on the dimension $n$, such that (see \cite{Beres-Lions})
\begin{equation} \label{stima-fuori}
|u_{k}(x)| \le C \, \frac{\| u_{k} \|_{H^{1}(\R^{n})}}{|x|^{\frac{n-1}{2}}} \qquad \forall\, |x|\ge \rho
\end{equation}
Moreover the hylomorphic solitons $\psi_k=u_k(x)e^{-i\omega_k t}$ satisfy the Dirichlet problems
\begin{equation} \label{problema-dentro}
\left\{
\begin{array}{ll}
-\Delta u_{k} = \omega_{k}^{2} u_{k} - W'(u_{k}) & \mbox{in } B(0,\rho)=\set{x\in \R^{n}\ :\ |x|\le \rho} \\[0.2cm]
u_{k}(x) = \eps_{k} &  \mbox{for } x\in \partial B(0,\rho)
\end{array} \right.
\end{equation}
where by (\ref{stima-fuori}), $(\eps_{k})$ are constants such that $\eps_{k} \to 0$ as $k\to \infty$. Now, recalling (\ref{NN}), without loss of generality (see Section 3 in \cite{BBBM}) we can assume that
\begin{equation} \label{crescita}
|N'(s)| \le c_{1} s^{p} + c_{2} s^{q}\ \mbox{ with }\ c_{1},c_{2}>0\ \mbox{ and }\ 1<p\le q < 2^{*}-1
\end{equation}
where $2^{*}= \frac{2n}{n-2}$. Hence, since $u_{k} \in H^{1}(B(0,\rho))$ and by the assumptions $\| u_{k} \|_{H^{1}} \to 0$ as $k\to \infty$, applying the classical Sobolev embedding theorems it follows that the right hand side of (\ref{problema-dentro}) satisfies
\begin{equation} \label{f-N}
f(u_{k}) := (\omega_{k}^{2} -1) u_{k} - N'(u_{k}) \longrightarrow 0 \quad \mbox{in } L^{s} \quad \forall\, s \in \left(1,\frac{2^{*}}{q} \right)
\end{equation}
where $q$ is given in (\ref{crescita}). Moreover, the classical theory of elliptic regularity (see for example \cite{adn}) implies that if the functions $u_{k}$ are solutions of (\ref{problema-dentro}) with $f(u_{k})$ in $L^{s}$ for some $s$, then
\begin{equation} \label{boot}
u_{k} \in W^{2,s} \quad \mbox{and} \quad \| u_{k} \|_{W^{2,s}} \le const \, \left( \| f_{k} \|_{L^{s}} + \eps_{k} \right)
\end{equation}
where the constant only depends on the radius $\rho$. Now, a classical bootstrap argument applies, and it follows that equations (\ref{f-N}) and (\ref{boot}) hold for bigger and bigger values of $s$. We sketch the ideas of the first step of the bootstrap. From (\ref{f-N}) and (\ref{boot}) it follows that
$$
\lim\limits_{k\to \infty}\, \| u_{k} \|_{W^{2,s}} =0 \quad \forall\, s \in \left(1,\frac{2^{*}}{q} \right)
$$
which implies that we can write (\ref{f-N}) for $s$ in an interval $(1,q_{1})$ with $q_{1}>\frac{2^{*}}{q}$, just applying again the classical Sobolev estimates. But then $f(u_{k})$ belongs to $L^{s}$ for all $s \in (1,q_{1})$, then (\ref{boot}) hold for all $s \in (1,q_{1})$. This argument can be repeated over and over until equations (\ref{f-N}) and (\ref{boot}) hold for all $s>1$.

In particular it follows that (\ref{f-N}) and (\ref{boot}) hold for $s>\frac n 2$, and the Sobolev theorems imply that
\begin{equation} \label{stima-dentro}
\| u_{k} \|_{L^{\infty}(B(0,\rho))} \le const \, \| u_{k} \|_{W^{2,s}(B(0,\rho))} \le const \, \left( \| f_{k} \|_{L^{s}(B(0,\rho))} + \eps_{k} \right)
\end{equation}
Hence, putting together (\ref{stima-fuori}) and (\ref{stima-dentro}) it follows that
$$
\lim\limits_{k\to \infty}\, \| u_{k} \|_{L^{\infty}(\R^{n})} = 0
$$
and the theorem is proved. \qed

\begin{definition}
Let $W$ satisfy (W-i), (W-ii) and (W-iii). We say that equation (NKG) (or W) is of
\emph{type ($\beta$)} if there exists $\beta_{0}>0$ such that any hylomorphic soliton
$\psi (t,x)$ of the form (\ref{sw}) satisfies
$$
\left\| \psi(t,\cdot )\right\|_{L^{\infty }(\R^{n})} \le \beta_{0} \qquad \forall\, t \in \R
$$
\end{definition}

\begin{theorem}
Assume that $W$ satisfies (W-i), (W-ii) and (W-iii), and write it as in (\ref{NN}). If
$$
N^{\prime }(s) \ge 0 \qquad \forall\, s \in \left(\beta _{0},+\infty \right)
$$
then equation NKG is of type ($\beta $).
\end{theorem}

\noindent {\bf Proof. } Let $\psi(t,x)=u(x)e^{-i\omega t}$ be a hylomorphic soliton, then $u$ is a solution of (\ref{static}). Setting $u=\beta_{0}+v$ it is sufficient to prove that $v(x)\le 0$. Let $\Omega :=\set{ x\ :\ v\left( x\right) >0}$. Substituting in (\ref{static}) we get
$$\left\{
\begin{array}{ll}
-\Delta v+W^{\prime }(\beta _{0}+v) = \omega^{2}(\beta_0+v) & \mbox{in} \ \Omega \\[0.2cm]
v = 0 & \mbox{in} \ \partial \Omega
\end{array}
\right.
$$
Multiplying both sides by $v$ and integrating in $\Omega $, we get
$$
0 = \int_{\Omega }\left[ \left| \nabla v\right| ^{2}+(1-\omega^2)(\beta_0+v)v+N^{\prime }(\beta_{0}+v)v-\omega^{2}v^{2}\right] dx \ge \int_{\Omega } (-\omega^{2}v^{2}) dx
$$
since $\omega^2<1$ by Theorem \ref{freq-adm}, and by assumption on $N'$. This implies that $v=0$ in $\Omega$. \qed

\begin{definition}
Let $W$ satisfy (W-i), (W-ii) and (W-iii). We say that equation (NKG) (or W) is of
\emph{type ($non$-$\beta$)} if for any $\beta _{0}>0$ there exists a hylomorphic soliton $\psi (t,x)$ of the form (\ref{sw}) for which
$$
\left\| \psi(t,\cdot )\right\|_{L^{\infty }(\R^{n})} \ge \beta_{0} \qquad \forall\, t \in \R
$$
\end{definition}

\begin{theorem}
Assume that $W$ satisfies (W-i), (W-ii) and (W-iii), and that
$$
\lim_{s\rightarrow \infty } \frac{W(s)}{s^2}=0 \quad \mbox{ and } \quad W(s)>0, \ \ s\not= 0
$$
Then equation (NKG) is of type ($non$-$\beta $).
\end{theorem}

\noindent {\bf Proof.} The argument is by contradiction. Let us assume that there exists a constant $C$ such that all hylomorphic solitons have $L^{\infty}$ norm bounded by $C$. Then we can choose a sequence of hylomorphic solitons $(u_k,\omega_k)$ with charges $\HH(u_{k},\omega_{k})$ diverging and with the properties
\begin{eqnarray}
&&\int u_k^2\, dx \rightarrow \infty \label{lim-2} \\
&& \| u_k(x) \|_{L^{\infty}(\R^{n})}\le C \label{lim-norm}
\end{eqnarray}
Let us now consider the sequence of radially symmetric triangle shaped functions $v_{k}$ defined by
\begin{equation} \label{trian}
v_{k} (r)= \left\{
\begin{array}{cl}
-\frac{M_k}{\rho_k}\, r + M_k & \mbox{ if $0 \leq r \le \rho_k$} \\[0.2cm]
0 & \mbox{ if $|r| \geq \rho_k$}
\end{array} \right.
\end{equation}
with parameters $M_{k}$ and $\rho_{k}$ suitable to satisfy $\int u_k^2dx=\int v_k^2dx$. For these functions it holds $\HH(v_{k},\omega_{k}) = \HH(u_{k},\omega_{k})$.

We now use (\ref{hyl}) and (\ref{alpha}) to obtain a contradiction. By the assumptions on $W$ and (\ref{lim-2}) and (\ref{lim-norm}), the sequence $(u_{k})$ satisfies
$$
\Lambda(u_{k},\omega_{k}) \ge \frac 1 2\, \omega_{k} + \frac 1 2 \, \frac{\int W(u_k)\, dx}{\int u_k^2\, dx} \succsim \frac 1 2\, \omega_{k} + \frac 1 2 \, \frac{W(C)}{C^{2}}
$$
letting the functions $u_{k}$ converging to the constant function $C$.

On the contrary, by (\ref{trian}) the functions $v_{k}$ satisfy
$$
\Lambda(v_{k},\omega_{k}) = \frac 1 2\, \omega_{k} + \frac 1 2 \, \frac{\int (|\nabla v_k| ^2 + W(v_k))\, dx}{\int v_k^2\, dx} = \frac 1 2\, \omega_{k} + O\left( \frac{W(M_{k})}{M_{k}^{2}} + {\frac{1}{\rho_{k}^2}}\right)
$$
Hence, letting $M_{k}\to \infty$ and $\rho_{k} \to \infty$ under the condition $\int u_k^2dx=\int v_k^2dx$, we get that
$$
\Lambda(v_{k},\omega_{k}) \le \Lambda(u_{k},\omega_{k})
$$
for $k$ sufficiently large. Hence, since $\HH(v_{k},\omega_{k}) = \HH(u_{k},\omega_{k})$,
$$
\E(v_{k},\omega_{k}) \le \E(u_{k},\omega_{k})
$$
for $k$ sufficiently large, which contradicts the fact that $(u_{k},\omega_{k})$ are points of absolute minimum for the energy. \qed

\bigskip

Given the definitions of ``types'', we classify the nonlinear Klein-Gordon equations in terms of the following classes

\begin{itemize}
\item \emph{($\alpha$, $\beta$)}; if it is of type ($\alpha $) and type ($\beta $)

\item \emph{($\alpha$, non-$\beta$)}; if it is of type ($\alpha $) and type \emph{(non-$\beta$)}

\item \emph{(non-$\alpha$, $\beta$)}; if it is type \emph{(non-$\alpha$)} and of type ($\beta $)

\item ($\gamma $); if it is of type \emph{(non-$\alpha$)} and of type \emph{(non-$\beta$)}
\end{itemize}

In particular the equation (NKG) is of type ($\gamma $) if for every $\gamma
_{0}>0$ there exist hylomorphic solitons $\psi _{s}(t,x)$ and $\psi _{b}(t,x)$
satisfying
$$
\left\| \psi_{s}(t,\cdot )\right\|_{L^{\infty }(\R^{n})} \le \gamma_{0}\leq \left\| \psi_{b}(t,\cdot) \right\|_{L^{\infty }(\R^{n})} \qquad \forall\, t \in \R
$$

We now give examples of nonlinear terms for all classes,
presenting numerical results obtained by means of the strategy discussed in Section \ref{minimizzazione_numerica}. The numerical investigation is carried out in a 2D context for simplicity.

\begin{figure}
    \begin{center}
    \subfigure[]
    {\includegraphics[width=6.0cm,height=4.5cm]{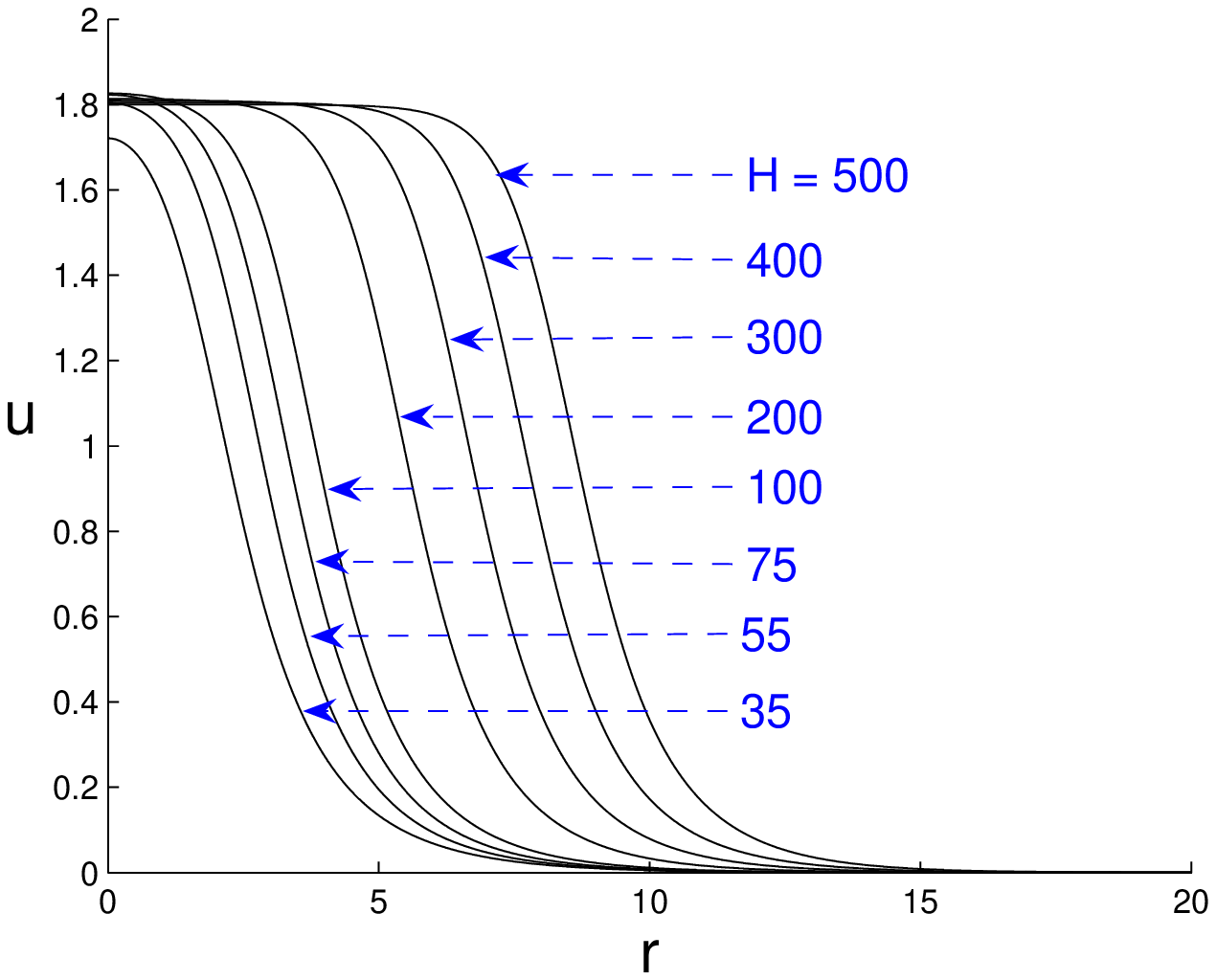}}
    \subfigure[]
    {\includegraphics[width=6.0cm,height=4.5cm]{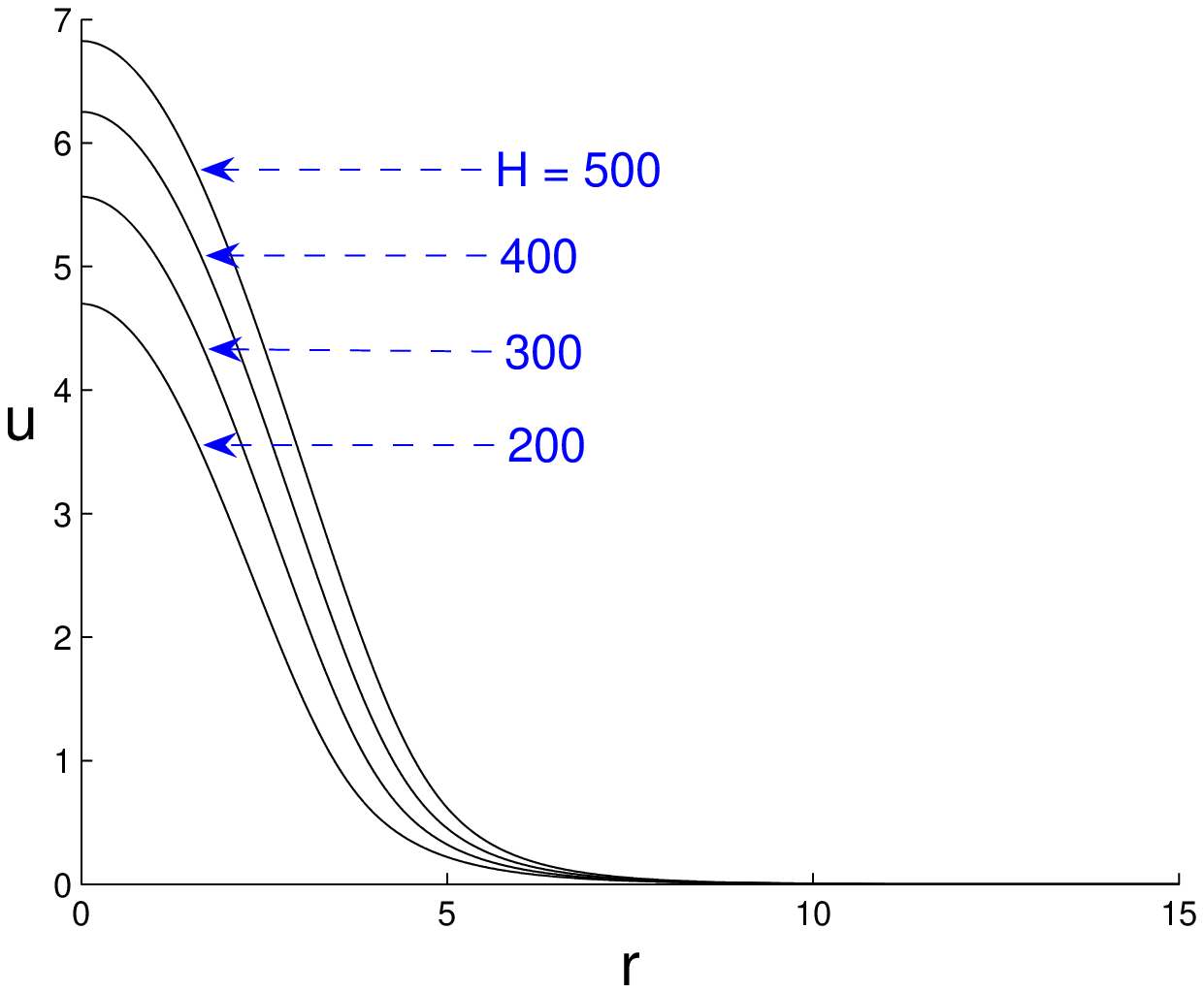}}
    \subfigure[]
    {\includegraphics[width=6.0cm,height=4.5cm]{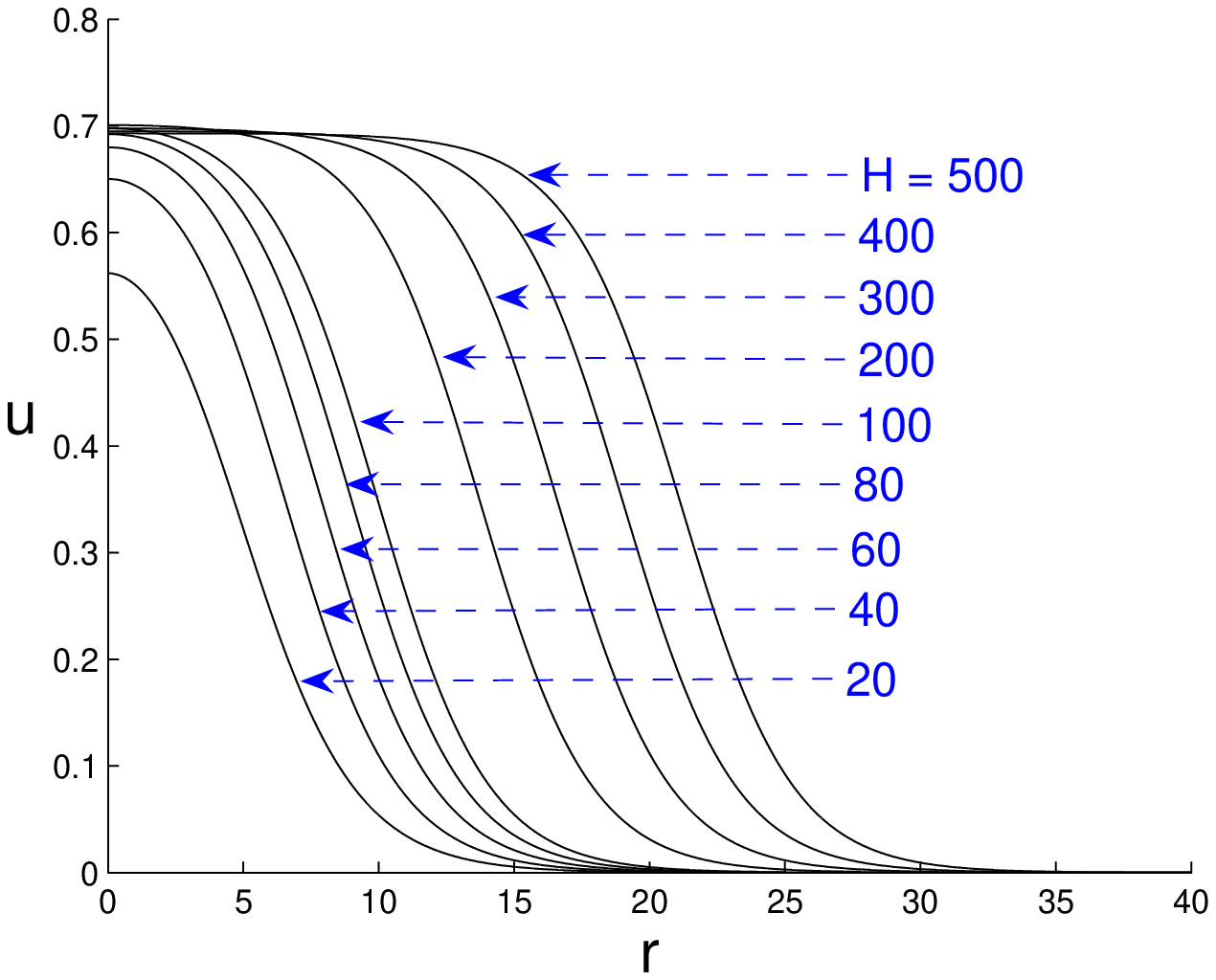}}
    \subfigure[]
    {\includegraphics[width=6.0cm,height=4.5cm]{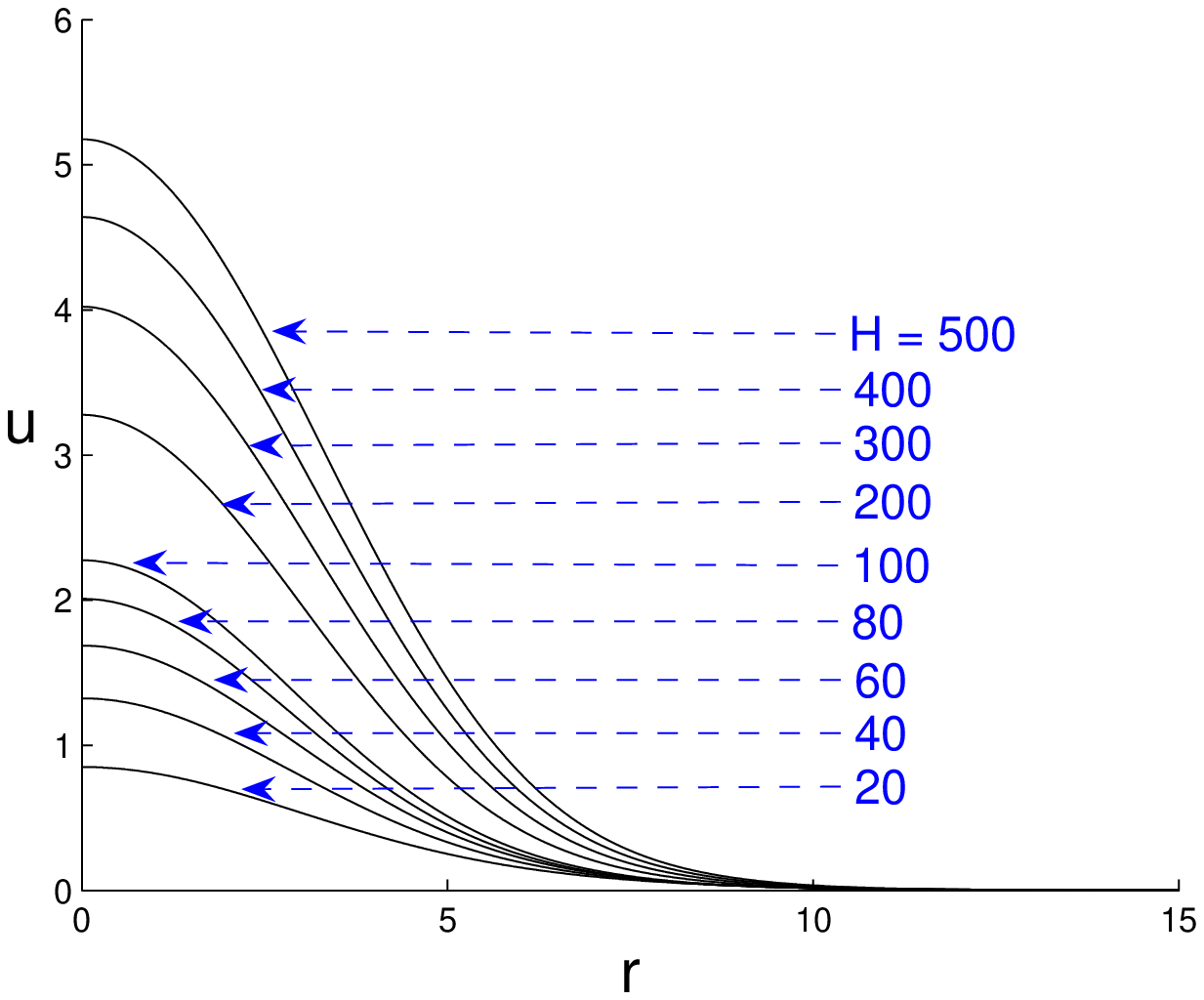}}
    \caption{Radial profile $u(r)$ of 2D Q-balls, for given values of the hylomorphic charge $\HH$.
             (a) W of type ($\alpha$, $\beta$), given by (\ref{ex_W_alfabeta}) with $a=2.5$;
           (b) W of type ($\alpha$, non-$\beta$), given by (\ref{ex_W_alfanonbeta});
           (c) W of type (non-$\alpha$,$\beta$), given by (\ref{ex_W_nonalfabeta}) with $a=1$;
           (d) W of type ($\gamma$), given by (\ref{ex_W_gamma}).}
    \label{fig_profili_radiali_u}
    \end{center}
\end{figure}

\begin{figure}
    \begin{center}
    \subfigure[]
    {\includegraphics[width=6.0cm,height=3.8cm]{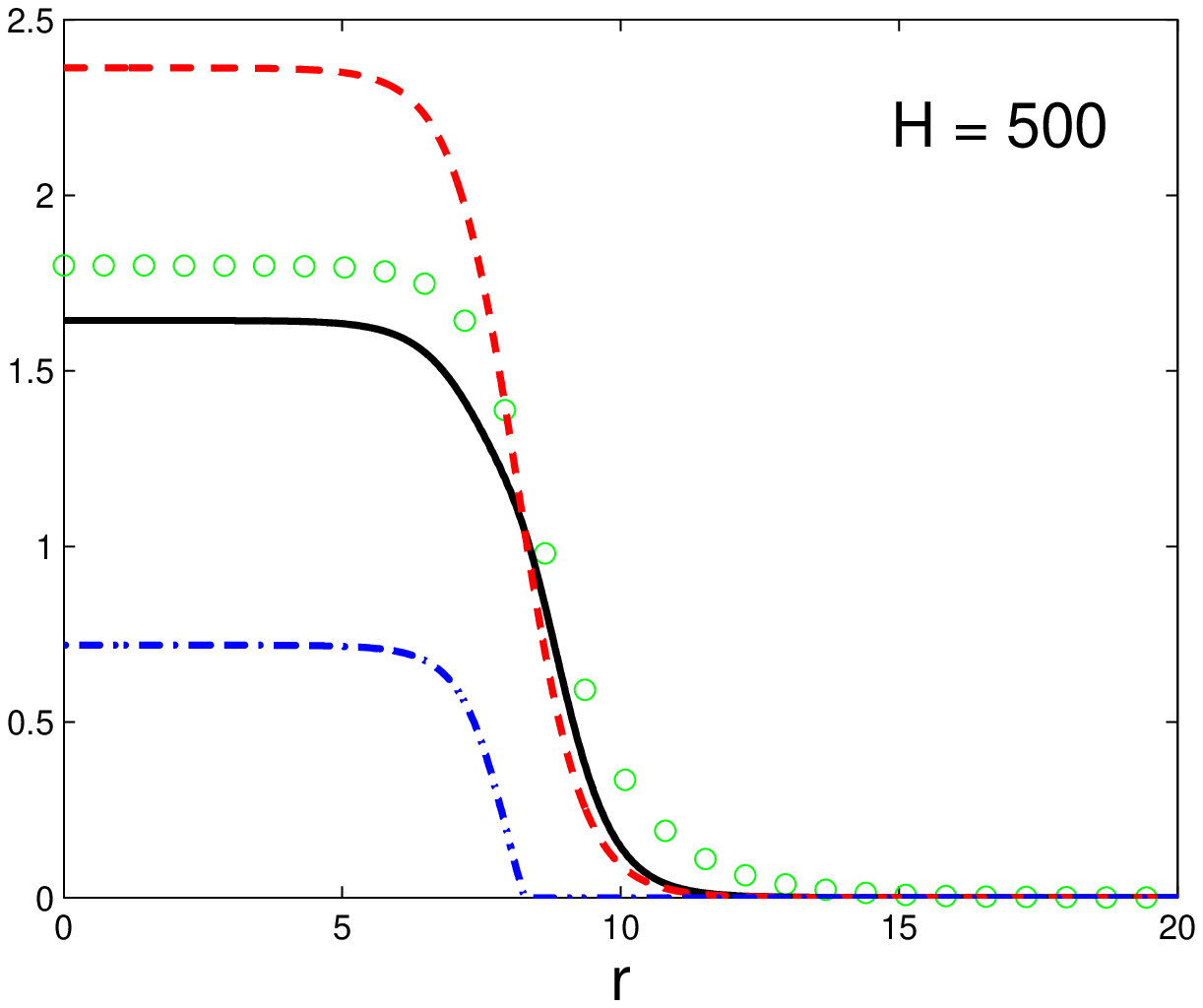}}
    \subfigure[]
    {\includegraphics[width=6.0cm,height=3.8cm]{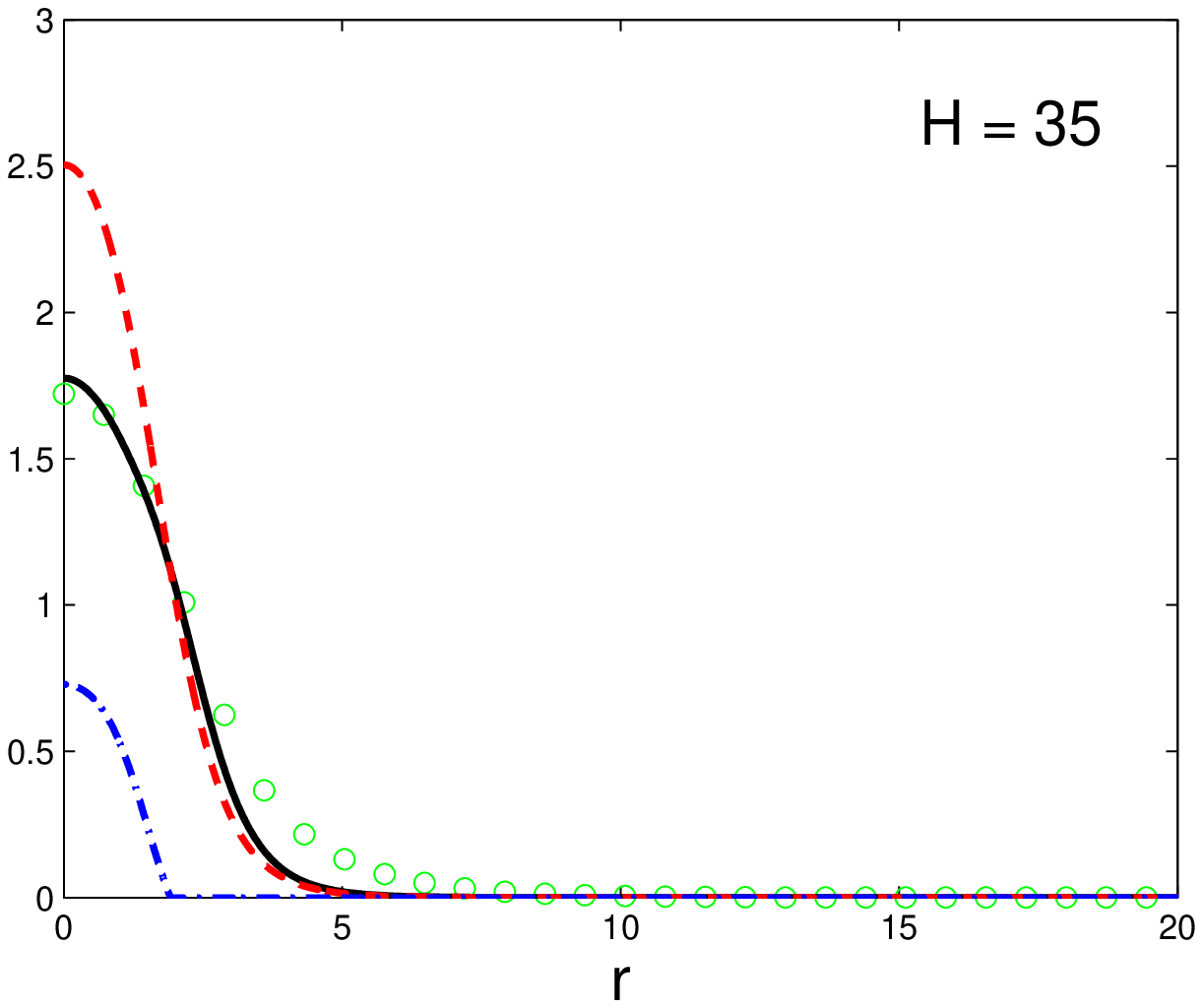}}
    \subfigure[]
    {\includegraphics[width=6.0cm,height=3.8cm]{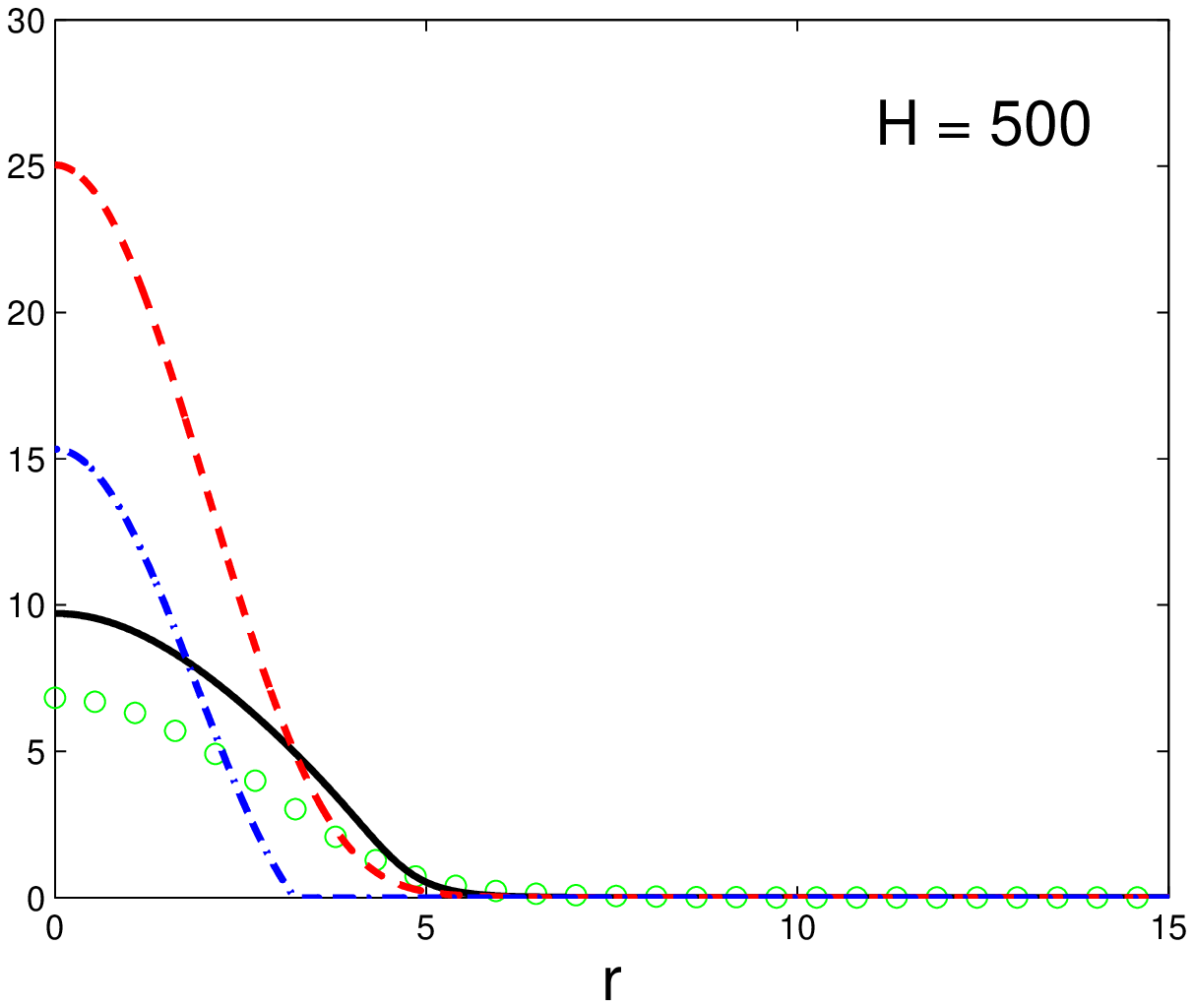}}
    \subfigure[]
    {\includegraphics[width=6.0cm,height=3.8cm]{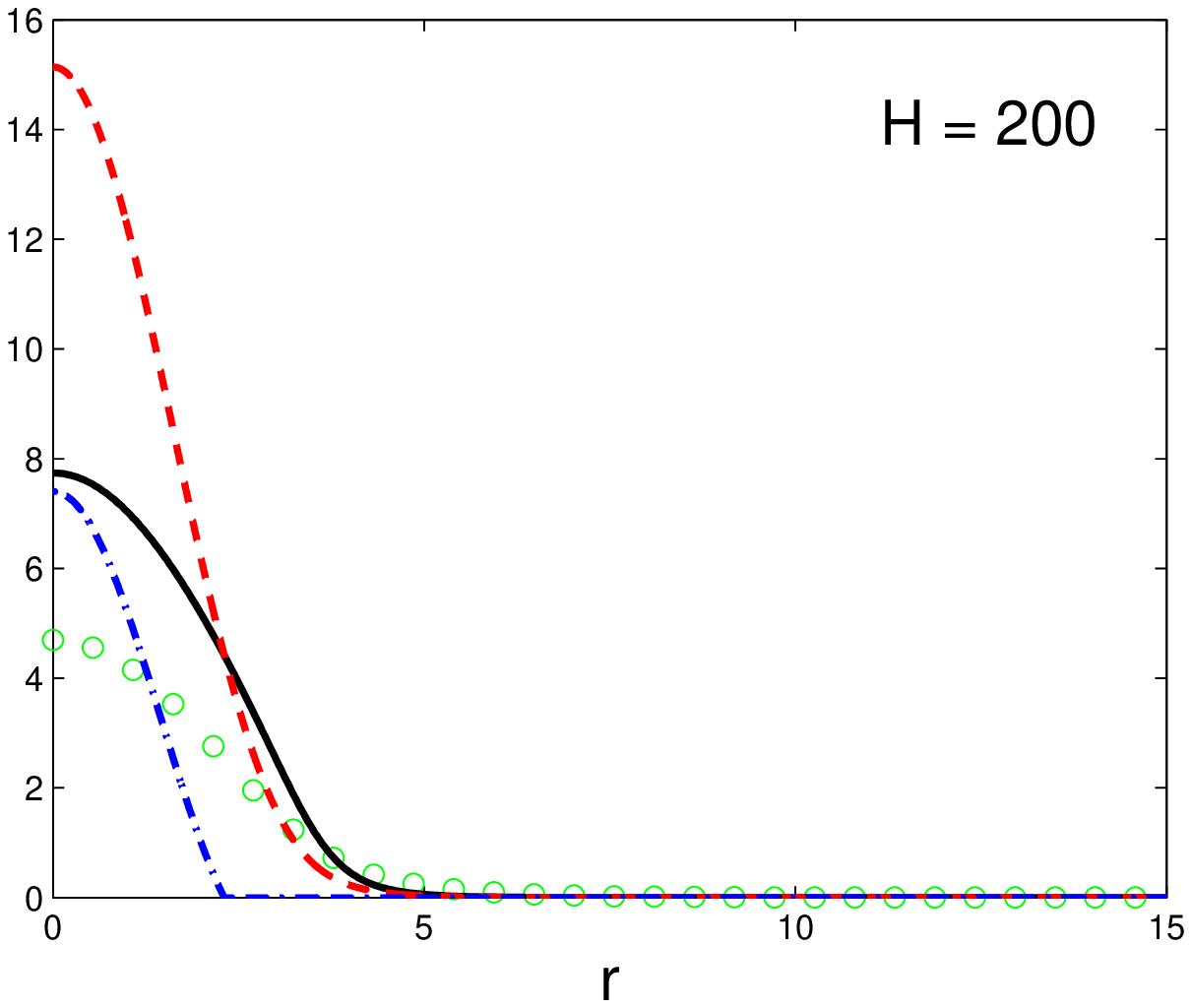}}
    \subfigure[]
    {\includegraphics[width=6.0cm,height=3.8cm]{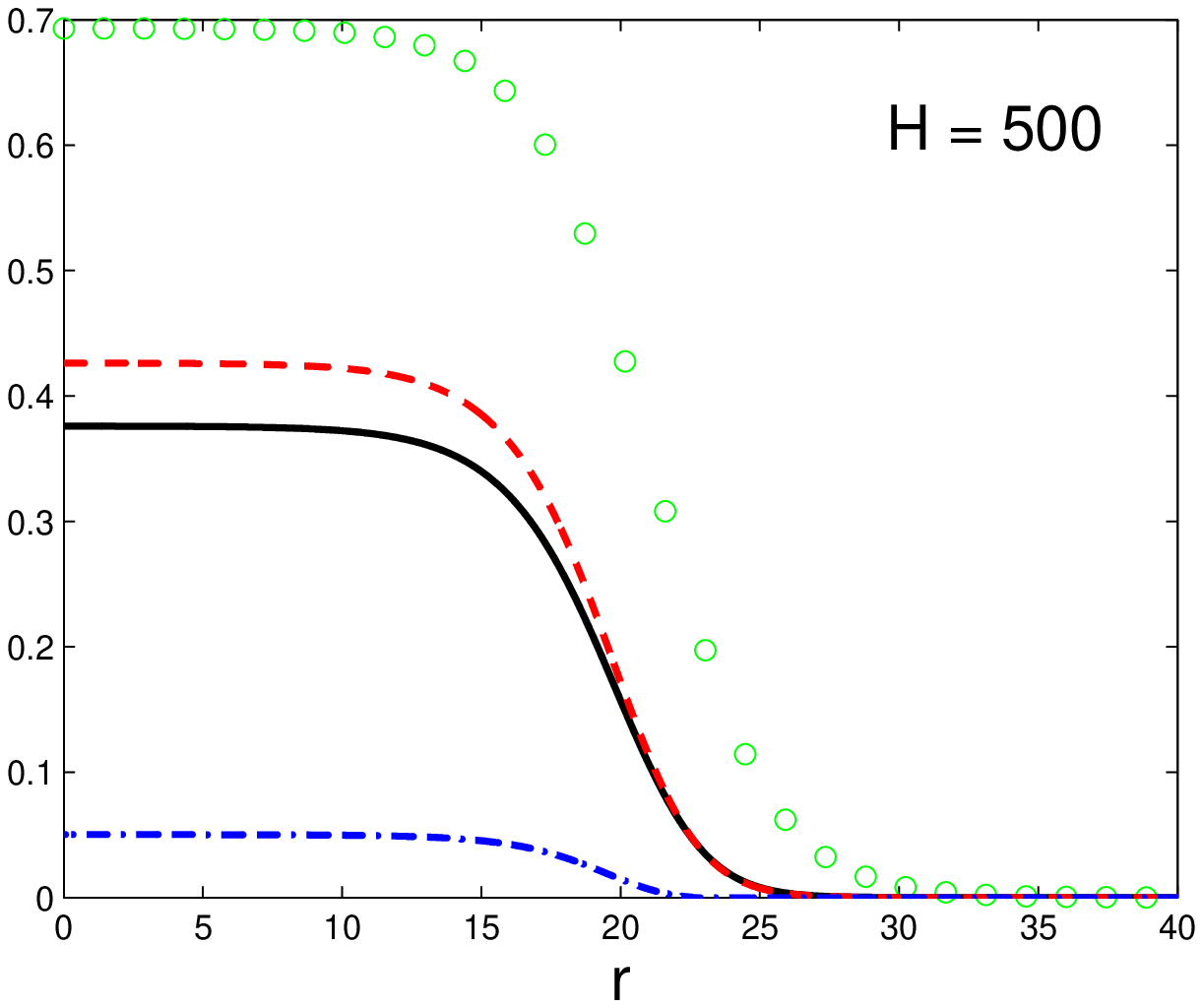}}
    \subfigure[]
    {\includegraphics[width=6.0cm,height=3.8cm]{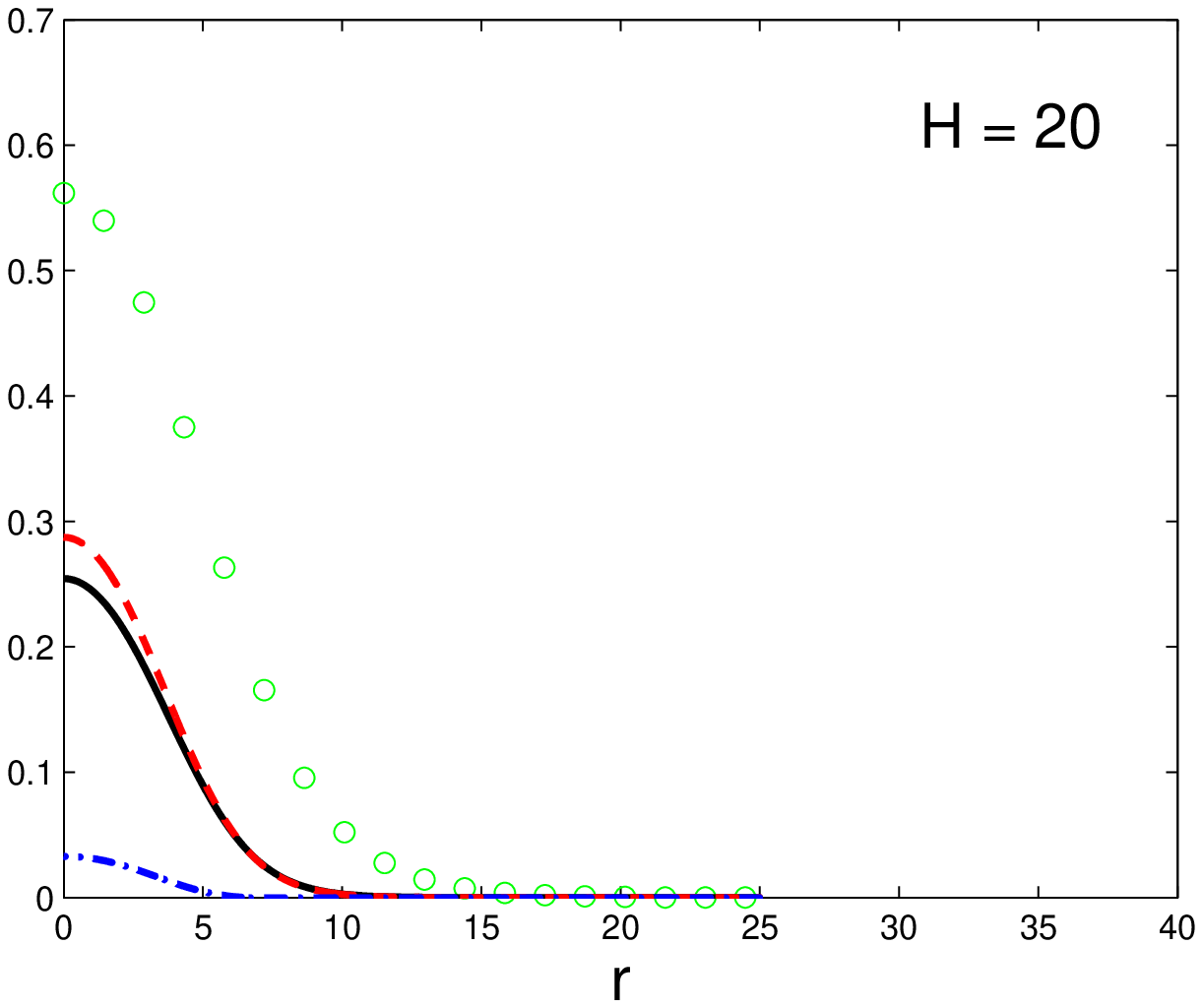}}
    \subfigure[]
    {\includegraphics[width=6.0cm,height=3.8cm]{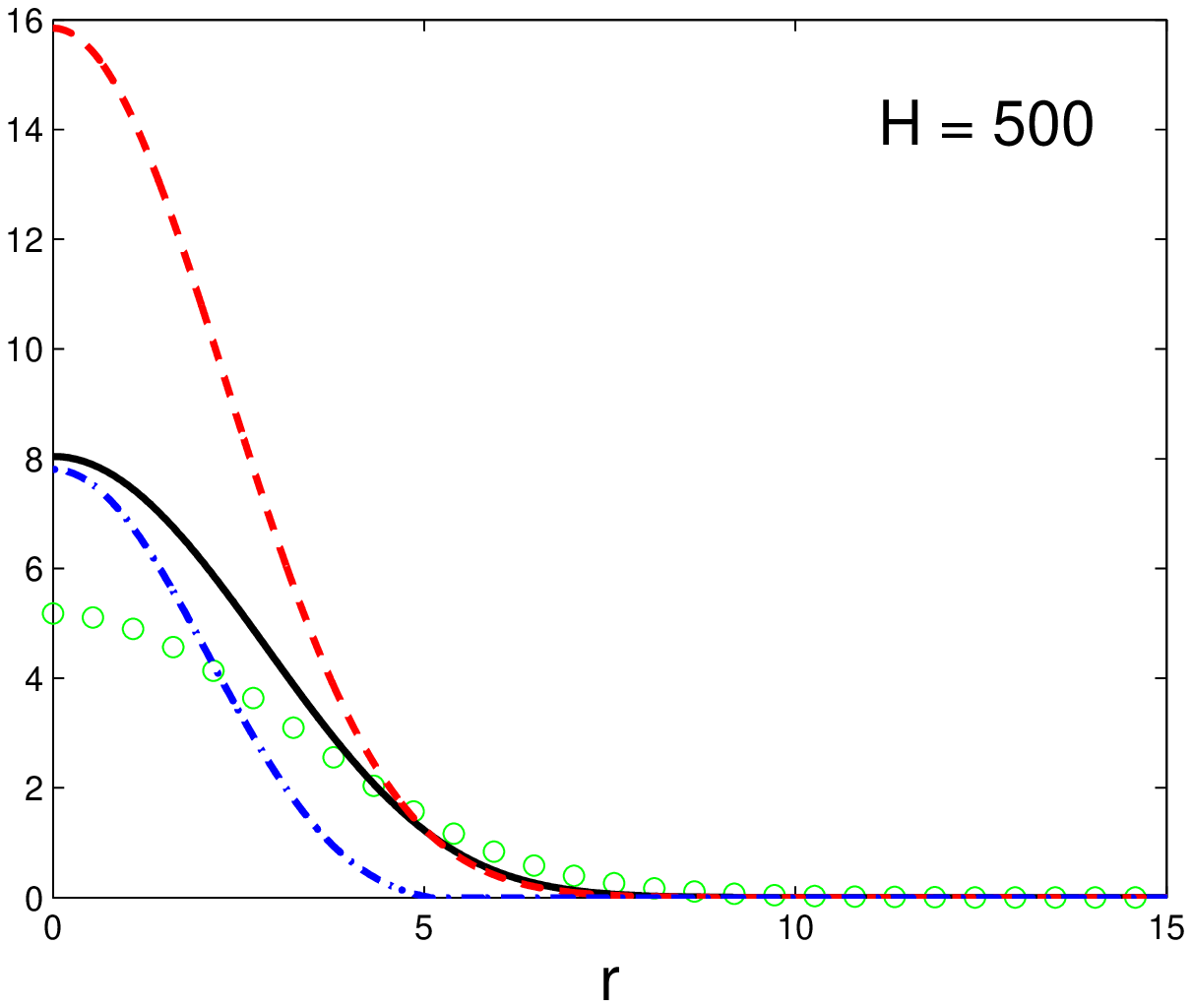}}
    \subfigure[]
    {\includegraphics[width=6.0cm,height=3.8cm]{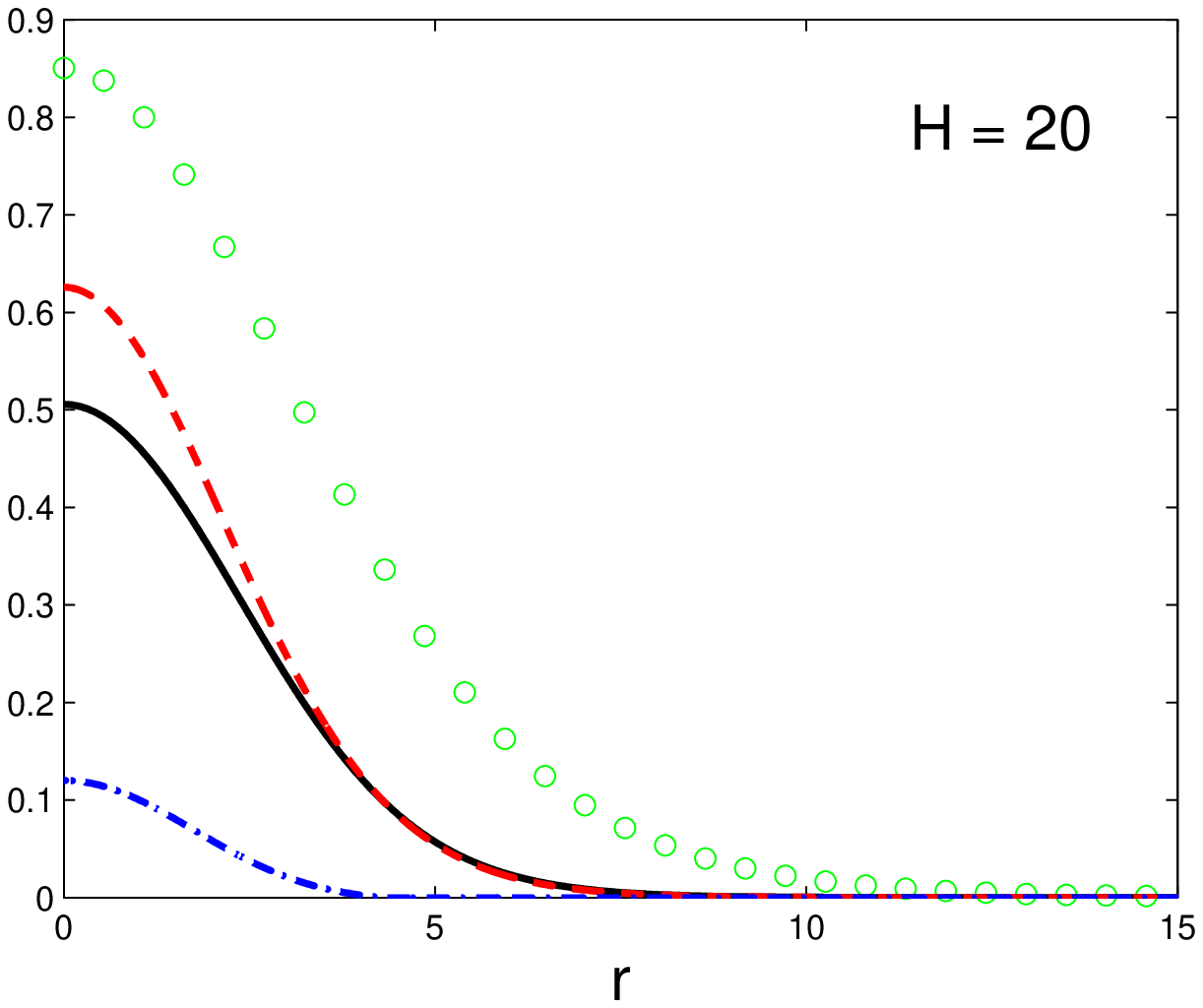}}
    \caption{Energy density ($\rho_{\E,\Psi}$, solid),
             hylomorphic charge density ($\rho_{\HH,\Psi }$, dashed) and
             corresponding binding energy density ($\rho_{B,\Psi}$, dashed-dotted)
           along the radial direction for 2D Q-balls having a given hylomorphic
           charge $\HH$ (specified in each figure).
             The radial profile is also shown through circles, for reference.
             (a) and (b) W of type ($\alpha$, $\beta$), given by (\ref{ex_W_alfabeta}) with $a=2.5$;
           (c) and (d) W of type ($\alpha$, non-$\beta$), given by (\ref{ex_W_alfanonbeta});
           (e) and (f) W of type (non-$\alpha$,$\beta$), given by (\ref{ex_W_nonalfabeta}) with $a=1$;
           (g) and (h) W of type ($\gamma$), given by (\ref{ex_W_gamma}).}
    \label{fig_supporti}
    \end{center}
\end{figure}

In particular, in Figure \ref{fig_profili_radiali_u} we have
plotted the radial profiles $u(r)$ of some $Q$-balls associated to
the (NKG) equations of the four classes. For each class, the
plotted profiles correspond to different values of hylomorphic
charge. We remark that for (NKG) equations of type ($\alpha$),
hylomorphic solitons with arbitrary small hylomorphic charge do
not exist. In particular, for the choice (\ref {ex_W_alfabeta}),
the threshold is approximately at $\sigma=25$. In this case, if
the charge is e.g. $\sigma=15$, the numerical algorithm leads to a
radial profile which is not a solution of the Klein-Gordon
equation in the unbounded domain, but rather a solution of the
Dirichlet problem in the interval fixed for the numerical
analysis. For what concerns hylomorphic solitons with big charges,
we point out the difference in the shape of the radial profiles
for (NKG) equations of type ($\beta$) or \emph{(non-$\beta$)}. In
the first case, as expected, there exists a bound for the
$L^{\infty}$ norm of the radial profiles, whereas in the second
case the $L^{\infty}$ norms increase arbitrarily. Finally in
Figure \ref{fig_supporti} we have plotted the energy density
($\rho_{\E,\Psi}$, solid), the hylomorphic charge density
($\rho_{\HH,\Psi }$, dashed) and the corresponding binding energy
density ($\rho_{B,\Psi}$, dashed-dotted) for some $Q$-balls of all
the four classes.

\bigskip

%--------------------------------------------------[ex1: alfa beta]------

\noindent
\textbf{Example no. 1: $W$ of type ($\alpha, \beta $).}
The equation
\begin{equation*}
\psi _{tt}-\Delta \psi +\left( 1+\left| \psi \right| -a\;\left| \psi \right|
^{2}+\left| \psi \right| ^{3}\right) \psi =0
\end{equation*}
is of type ($\alpha, \beta $) provided that $a\cong 2.5.$
Here
\begin{equation}
W \left( \psi \right) =
\frac{1}{2}\left| \psi \right| ^{2}+\frac{1}{3}\left|
\psi \right| ^{3}-\frac{a}{4}\;\left| \psi \right| ^{4}+\frac{1}{5}\left|
\psi \right| ^{5}
\label{ex_W_alfabeta}
\end{equation}
We have $\alpha_0\cong 1$ and $\beta_0\cong 2.5.$

%------------------------------------------[ex2: alfa non-beta]--------------

\bigskip
\noindent
\textbf{Example no. 2: $W$ of type \emph{($\alpha$, non-$\beta$)}.}
The equation
$$
\psi _{tt}-\Delta \psi +\frac{\psi }{1-\left| \psi \right| +\left| \psi \right|^{2}}=0
$$
is of type \emph{($\alpha$, non-$\beta$)}. Here
\begin{equation}
W \left( \psi \right) =
\frac{1}{2}\log \left( 1-\left| \psi \right| +\left|
\psi \right| ^{2}\right) +\frac{1}{\sqrt{3}}\left[ \arctan \left( \frac{%
2\left| \psi \right| -1}{\sqrt{3}}\right) +\frac{\pi }{6}\right]
\label{ex_W_alfanonbeta}
\end{equation}

%----------------------------------------------[ex.3 non-alfa beta]---------

\bigskip
\noindent
\textbf{Example no. 3: $W$ of type \emph{(non-$\alpha$, $\beta$)}.}
The equation
$$
\psi _{tt}-\Delta \psi +\left( 1-a\;\left| \psi \right| +\left| \psi \right| ^{2}\right) \psi =0
$$
is of type \emph{(non-$\alpha$, $\beta$)} provided that $a\in \left( 0,2\right)$.
Here
\begin{equation}
W \left( \psi \right) =
\frac{1}{2}\left| \psi \right| ^{2}-\frac{a}{3}\left|
\psi \right| ^{3}+\frac{1}{4}\;\left| \psi \right| ^{4}
\label{ex_W_nonalfabeta}
\end{equation}

%--------------------------------------------[ex.4 gamma] ---------------

\bigskip
\noindent
\textbf{Example no. 4: $W$ of type ($\gamma$).}
The equation
$$
\psi _{tt}-\Delta \psi +\frac{\psi }{1+\left| \psi \right| }=0
$$
is of type ($\gamma $). Here
\begin{equation}
W \left( \psi \right) =
\left| \psi \right| -\log \left( 1+\left| \psi \right|
\right) = \frac 1 2\, |\psi|^{2} - \frac 1 3\, |\psi|^{3} + o(|\psi|^{3})
\label{ex_W_gamma}
\end{equation}
We notice that $W(s)$ fulfils the assumption (W-iv). The existence of hylomorphic solitons with arbitrary charge is guaranteed by Corollary \ref{123-carica-tutte}.

%-------------------------------------------[FIGURE]------------------------

\section{Acknowledgments}
The research of the fourth author has been partially funded by the Mathematics Dept. of the University of Bari (grant as for D.D. n.1 - 19/09/2007), which is most gratefully acknowledged.
%%%%%%%%%%%%%%%%%%%%%%%%%%%%%%%%%%%%%%%%%%%%%%%%%%%%%%%%%%%%%%%%%%%%%%%%%%%%%

\end{document}